\documentclass[12pt,oneside,reqno]{amsart}

\usepackage{amssymb}
\usepackage{algpseudocode}
\usepackage{algorithm}
\usepackage{verbatim}
\usepackage{cite}
\usepackage{sidecap}
\usepackage{wrapfig}
\usepackage{subfig}
\RequirePackage{dsfont} 
 \scrollmode

\usepackage{graphicx}

\usepackage[usenames]{color}
\definecolor{darkred}{RGB}{139,0,0}
\definecolor{darkgreen}{RGB}{0,100,0}
\definecolor{darkmagenta}{RGB}{139,0,139}

\def\rho{\varrho}

\def\rd{\,{\mathrm d}}

\theoremstyle{plain}
\newtheorem{theorem}{Theorem}
\newtheorem{lemma}{Lemma}
\newtheorem{corollary}{Corollary}
\newtheorem{proposition}{Proposition}

\theoremstyle{definition}
\newtheorem{remark}{Remark}

\begin{document}

\title
[Optimal approximation of stochastic integrals in analytic  noise model]
{Optimal approximation of stochastic integrals in analytic noise model}


\author[A. Kaluza]{Andrzej Ka\l u\.za}
\address{AGH University of Science and Technology,
Faculty of Applied Mathematics,
Al. A.~Mickiewicza 30, 30-059 Krak\'ow, Poland}
\email{akaluza@agh.edu.pl}

\author[P. Morkisz]{Pawe{\l} M. Morkisz}
\address{AGH University of Science and Technology,
Faculty of Applied Mathematics,
Al. A.~Mickiewicza 30, 30-059 Krak\'ow, Poland}
\email{morkiszp@agh.edu.pl}

\author[P. Przyby{\l}owicz]{Pawe{\l} Przyby{\l}owicz}
\address{AGH University of Science and Technology,
Faculty of Applied Mathematics,
Al. A.~Mickiewicza 30, 30-059 Krak\'ow, Poland}
\email{pprzybyl@agh.edu.pl, corresponding author}

\begin{abstract}
We study  approximate stochastic It\^o integration of processes belonging to a class of progressively measurable stochastic processes that are H\"older continuous in the $r$th mean. 

Inspired by increasingly popularity of computations with low precision (used on Graphics Processing Units -- GPUs and standard Computer Processing Units -- CPU for significant speedup), we introduce a suitable analytic noise model of standard noisy information about $X$ and $W$. In this model we show that the upper bounds on the  error of the Riemann-Maruyama quadrature are proportional to $n^{-\varrho}+\delta_1+\delta_2$, where $n$ is a number of noisy evaluations of $X$ and $W$, $\varrho\in (0,1]$ is a H\"older exponent of $X$, and $\delta_1,\delta_2\geq 0$ are precision parameters for values of $X$ and $W$, respectively. Moreover, we show that the error of any algorithm based on at most $n$ noisy evaluations of $X$ and $W$ is at least  $C(n^{-\varrho}+\delta_1)$.  Finally, we report numerical experiments performed on both CPU and GPU, that confirm our theoretical findings,  together with some computational performance comparison between those two architectures.
\newline
\newline
\textbf{Key words:} Wiener process, noisy information,  analytic noise model, optimal approximation, minimal error, GPU
\newline
\newline
\textbf{Mathematics Subject Classification:} 68Q25, \ 65C30.

\end{abstract}
\maketitle

\section{Introduction}\label{sec:intr}
In this paper we investigate the problem of optimal approximation of stochastic integrals of the following form
\begin{equation}
	\label{PROBLEM}
		\mathcal{I}(X,W)=\int\limits_0^T X(t)\rd W(t),
\end{equation}
where $T>0$ and $W=\{W(t)\}_{t\geq 0}$ is a one-dimensional Wiener process on some probability space 
$(\Omega,\Sigma,\mathbb{P})$, and we consider integrands $X=\{X(t)\}_{t\in [0,T]}$ from  a class of progressively measurable stochastic processes that are H\"older continuous in the $r$th mean. Such quadrature problems arise, for example, in the context of numerical solutions of stochastic differential equations, see Section 4.4 in \cite{KP} and \cite{esikr}. Since the exact values of such stochastic integrals are known only in limited cases, an efficient approximation of $\mathcal{I}(X,W)$ with the  error as small as possible is of interest. We are aiming at methods that are based only on discrete values of $X$ and $W$ which are, additionally, corrupted with some noise. 

The problem of optimal approximation of stochastic It\^o integrals under exact information about $X$ and $W$ has been well studied in the literature, see, \cite{hein1}, \cite{heindaun}, \cite{hert}, \cite{PP1}, \cite{PP2}, \cite{PP3}, \cite{waswo}. Less explored is the approximation of stochastic integrals in the case when  values of $X$ and $W$  are  corrupted with some noise. The noise may arise from measurement errors, previous computations or simply the floating point numbers representation errors. There is a trend observed in deep learning, where the computations are mostly conducted in lower precision, i.e. using not only single but also half precision. Exemplary performance for NVIDIA V100 graphic card is about 7 TFLOPS for double precision and 14 TFLOPS for single precision. For deep learning purposes, using the nature of the operations and also further lowering the precision to half precision in chosen operations, enabled obtaining up to 112 TFLOPS of performance \cite{NV_Volta}. Hence, it is a huge motivation for analyzing how lowering the precision for other than deep learning applications, influence the computed result.

There are many results on solving problems under noisy information, including the problems such as integrating or approximating  regular functions (see, e.g., \cite{JoC4, Pla96, KaPl90}), $L_p$ approximation of piecewise regular functions (see \cite{MoPl16}), approximate solving of IVPs (see \cite{KaPr16}) or PDEs (see \cite{Wer96, Wer97}). 

In this paper we study noisy information for stochastic It\^o integration. According to our best knowledge this is the first paper that deals with noisy information for stochastic processes and its application to numerical computation of It\^o integrals. In this sense this paper can be seen as the extension of the model presented in \cite{MoPr17} in the context of SDEs. However, in that paper only drift and diffusion coefficients were corrupted with noise, information about the Wiener process was exact.

We use the Information--Based Complexity framework (see \cite{Pla96} and \cite{TWW88}). We assume that the algorithms may use only  noisy standard information about the integrand $X$ and the Wiener process $W$. Namely, let  $\delta_1, \delta_2 \geq 0$ be the precision levels corresponding to the processes $X$ and $W$, respectively. (The case of $\delta_1 = \delta_2 = 0$ corresponds to the exact information.) Available standard information about each coefficient consists of  noisy evaluations of $X$ and $W$ of points $t_i, z_i \in [0,T]$. This means that, for example, for the Wiener process and for a given $z_i\in[0,T]$ an evaluation returns  $\tilde W(z_i)$ such that $| W(z_i) - \tilde W(z_i)| \leq \delta_2(1+|W(z_i)|^s)$ for some $s\geq 0$. In the context of computations performed on GPUs, this can be interpreted as the standard relative error. (Detailed description of the noisy  information is given in Section 2.) From the reasons that become clear in Section 2 we refer to the model presented as to {\it analytical noise model}.

The error of an algorithm is measured in the $r$-th mean  maximized over the class of input data $X$ and over all permissible information about $(X,W)$ with the given precisions $\delta_1, \delta_2\geq 0$. In the model we show that the upper bounds on the  error of the Riemann-Maruyama quadrature are proportional to $n^{-\varrho}+\delta_1+\delta_2$, where $n$ is a number of noisy evaluations of $X$ and $W$, $\varrho\in (0,1]$ is a H\"older exponent of $X$, and $\delta_1,\delta_2\geq 0$ are precision parameters for values of $X$ and $W$, respectively. Moreover, we show that the error of any algorithm based on at most $n$ noisy evaluations of $X$ and $W$ is at least  $C(n^{-\varrho}+\delta_1)$. We also present numerical experiments that confirm our theoretical findings. As we perform the similar operations for multiple trajectories of considered processes, the proposed algorithm is highly parallel. Hence, the usage of GPU for computation acceleration is very natural. There is a vast list of problems where employing GPUs gives significant speedups \cite{Ryoo}, including matrix operations \cite{Fatahalian, Kruger}, bioinformatics \cite{Langdon} or solving ordinary or random differential equations \cite{Riesinger, Riesinger2}.

The paper is organized as follows. Section 2 contains basic notion and definitions. The Riemann-Maruyama quadrature rule $\mathcal{A}_n^{RM}$ for the perturbed information and upper estimate on its error (Theorem 1) are presented in Section 3. Section 4 consists of some lower bounds (Lemma \ref{lower_bounds}, Proposition \ref{prop_lower}). This leads to the conclusion that the algorithm $\mathcal{A}_n^{RM}$ is optimal (Theorem \ref{OPT_RME_ST}). Section 5 reports
numerical experiments performed for the GPU implementation of the algorithm.
\section{Preliminaries}\label{sec:prel}
\noindent
Denote $\mathbb{N}=\{1,2,\ldots\}$. Let  $W = \{W(t)\}_{t\geq 0}$ be the standard, one-dimensional Wiener process defined on a complete probability space  $(\Omega,\Sigma, \mathbb{P})$. By $\{\Sigma_t\}_{t\geq 0}$ we denote a filtration, satisfying the usual conditions, such that $W$ is a Wiener process with respect to $\{\Sigma_t\}_{t\geq 0}$.

For a random variable $Y:\Omega\to\mathbb{R}$ we write $\|Y\|_q=(\mathbb{E}|Y|^q)^{1/q}$, $q\geq 2$. Moreover, by $\mathcal{L}$ we mean the following differential operator
\begin{equation}
	\mathcal{L}=\frac{\partial}{\partial t}+\frac{1}{2}\frac{\partial^2}{\partial y^2}.
\end{equation}
For $r,q\in [2,+\infty)$, $q\geq r$, $L\geq 0$, $\varrho\in (0,1]$ we consider the following class $F^{\varrho,r,q}_L$ of stochastic processes $X=\{X(t)\}_{t\in [0,T]}$ 

\begin{align*}
	F^{\varrho,r,q}_L=\{ X:[0,T]\times \Omega\to\mathbb{R} \ | \ &X \ \hbox{is} \ \{\Sigma_t\}_{t\geq 0}-\hbox{progressively measurable},\\
		& \Bigl\|\sup\limits_{t\in [0,T]} |X(t)|\Bigl\|_q\leq L,\\
		& \|X(t)-X(z)\|_r\leq L|t-z|^{\varrho}, \ t,z\in[0,T]\}.
\end{align*}

Let us recall that by Theorem 33 from Chapter IV in \cite{delmey} for $X$ being $\{\Sigma_t\}_{t\geq 0}$-progressively measurable, the process $\{\sup\limits_{0\leq z\leq t}|X(z)|\}_{t\in [0,T]}$ is also $\{\Sigma_t\}_{t\geq 0}$-progressively measurable. Hence, $\sup\limits_{t\in [0,T]} |X(t)|$
is ($\Sigma_T$-measurable) random variable. Moreover, the processes from the class $F^{\varrho,r,q}_L$ are It\^o integrable, see, for example, \cite{legall} and \cite{KARSHR}. 

The numbers $r,q,L,\varrho,T$   will be called parameters of the class $F^{\varrho,r,q}_L$. Except for $T$ the parameters are not known
and the algorithm presented later on will not use them as input parameters.

In order to define suitable model of computation under inexact information about $X$, $W$ we need to introduce the following auxiliary classes.
\newline
Let
\begin{align*}
	\mathcal{K}^1= & \, \{p:[0,T]\times\mathbb{R}\to\mathbb{R} \ | \ p(\cdot,\cdot)-\hbox{Borel measurable},\\ 
&  \; |p(t,y)|\leq 1+|y| \ \hbox{for all} \ t\in [0,T], y\in\mathbb{R}\}.
\end{align*}

For $s\in [0,+\infty)$  we define
\begin{align*}	
	\mathcal{K}^{2}_{s}=& \,\Bigl\{p:[0,T]\times\mathbb{R}\to\mathbb{R} \ | \ p\in C^{1,2}([0,T]\times\mathbb{R}),\\
& \;\max\Bigl\{\Bigl|\frac{\partial p}{\partial t}(t,y)\Bigl|,\Bigl|\frac{\partial p}{\partial y}(t,y)\Bigl|,\Bigl|\frac{\partial^2 p}{\partial y^2}(t,y)\Bigl|\Bigr\}\leq 1+|y|^s \ \hbox{for all} \ t\in [0,T], y\in\mathbb{R} \Bigr\},
\end{align*}
\begin{align*}	
\mathcal{\bar K}^{2}_{s}= & \, \Bigl\{p:[0,T]\times\mathbb{R}\to\mathbb{R} \ | \ p\in C^{1,1}([0,T]\times\mathbb{R}),\\
&\;\max\Bigl\{\Bigl|\frac{\partial p}{\partial t}(t,y)\Bigl|,\Bigl|\frac{\partial p}{\partial y}(t,y)\Bigl|\Bigr\}\leq 1+|y|^s \ \hbox{for all} \ t\in [0,T], y\in\mathbb{R} \Bigr\},
\end{align*}
(if $s=0$ then we set $|y|^s:=0$ for all $y\in\mathbb{R}$), 
and for $\alpha,\beta\in (0,1]$ 
\begin{align*}
\mathcal{K}^{3}_{\alpha,\beta}=&\,\{p:[0,T]\times\mathbb{R}\to\mathbb{R} \ | \ |p(t,x)-p(z,y)|\leq |t-z|^{\alpha}+|x-y|^{\beta}\notag\\
&\; \ \hbox{for all} \ t,z\in [0,T], x,y\in\mathbb{R} \}
\end{align*}
We have that $\mathcal{K}^2_0\subset\mathcal{K}^3_{1,1}$. In the sequel, the classes above will allow us to to model, at least in some sense, the influence of the regularity of the noise on the error bound.

For $\delta_1,\delta_2\geq 0$ we define
\begin{equation}
	V_X(\delta_1)=\{\tilde X \ | \ \exists_{p_X\in\mathcal{K}^1}:\forall_{(t,\omega)\in [0,T]\times\Omega} \ \tilde X(t,\omega) = X(t,\omega)+\delta_1\cdot p_X(t,X(t,\omega))\},
\end{equation}
and the classes of disturbed Wiener process
\begin{equation}
	\mathcal{W}_{s}(\delta_2)=\{\tilde W \ | \ \exists_{p_W\in\mathcal{K}^{2}_{s}}:\forall_{(t,\omega)\in [0,T]\times\Omega} \ \tilde W(t,\omega) = W(t,\omega)+\delta_2\cdot p_W(t,W(t,\omega))\},
\end{equation}
\begin{equation}
	\mathcal{\bar W}_{s}(\delta_2)=\{\tilde W \ | \ \exists_{p_W\in\mathcal{\bar K}^{2}_{s}}:\forall_{(t,\omega)\in [0,T]\times\Omega} \ \tilde W(t,\omega) = W(t,\omega)+\delta_2\cdot p_W(t,W(t,\omega))\},
\end{equation}
\begin{equation}
	\mathcal{W}_{\alpha,\beta}(\delta_2)=\{\tilde W \ | \ \exists_{p_W\in\mathcal{K}^{3}_{\alpha,\beta}}:\forall_{(t,\omega)\in [0,T]\times\Omega} \ \tilde W(t,\omega) = W(t,\omega)+\delta_2\cdot p_W(t,W(t,\omega))\}.
\end{equation}
Since we impose some functional structure on corrupting functions, we refer to this model as to {\it analytic noise model}. (See also Remark \ref{rem_1} for possible alternative approach.)

For $X\in F^{\varrho,r,q}_L$ let $\tilde X\in V_X(\delta_1)$ and $\tilde W \in\mathcal{W}(\delta_2)$ where  $\mathcal{W}\in\{\mathcal{W}_{s},\mathcal{W}_{\alpha,\beta}\}$. We assume that the algorithm is based on discrete noisy information about $X$ and $W$. Hence, a vector of noisy information has the following form
\begin{equation}
	\mathcal{N}(\tilde X,\tilde W)=[\tilde X(t_0),\tilde X(t_1),\ldots,\tilde X(t_{i_1-1}),\tilde W(z_0),\tilde W(z_1),\ldots,\tilde W(z_{i_2-1})],
\end{equation}
where $i_1,i_2\in\mathbb{N}$. Moreover, $t_0,t_1,\ldots,t_{i_1-1}\in [0,T]$ and $z_0, z_1,\ldots, z_{i_2-1}\in [0, T ]$ are
given time points. Hence, the information is nonadaptive (see \cite{hein1} and \cite{TWW88} for more discussion on adaptive and nonadaptive information). We assume that $t_i \neq t_j, z_i \neq z_j$ for all $i \neq j$. The total number of (noisy) evaluations of $X$ and $W$ is $l=i_1+i_2$.

An algorithm $\mathcal{A}$ using $\mathcal{N}(\tilde X,\tilde W)$, that approximates $\mathcal{I}(X,W)$, is of the form
\begin{equation}	
	\label{appr_met}
	\mathcal{A}(\tilde X,\tilde W,\delta_1,\delta_2)=\varphi(\mathcal{N}(\tilde X,\tilde W)),
\end{equation}
where
\begin{equation}
	\varphi:\mathbb{R}^{i_1+i_2}\to\mathbb{R},
\end{equation}
is a Borel measurable mapping.

For a given $n\in\mathbb{N}$ we denote by $\Phi_n$ a class of all algorithms of the form (\ref{appr_met}) for which the total number of evaluations $l$
is at most $n$.

For a fixed $X\in F^{\varrho,r,q}_L$ the error of $\mathcal{A}\in\Phi_n$ is defined as
\begin{equation}
	e^{(r)}(\mathcal{A},X,\mathcal{W},\delta_1,\delta_2)=\sup\limits_{(\tilde X,\tilde W)\in V_X(\delta_1)\times \mathcal{W}(\delta_2)}\|\mathcal{I}(X,W)-\mathcal{A}(\tilde X,\tilde W,\delta_1,\delta_2)\|_r,
\end{equation} 
where $\mathcal{W}\in\{\mathcal{W}_{s},\mathcal{\bar W}_{s},\mathcal{W}_{\alpha,\beta}\}$.
The worst case error of $\mathcal{A}$ in $\mathcal{G}$ is given by
\begin{equation}
	e^{(r)}(\mathcal{A},\mathcal{G},\mathcal{W},\delta_1,\delta_2)=\sup\limits_{X\in \mathcal{G}}e^{(r)}(\mathcal{A},X,\mathcal{W},\delta_1,\delta_2)
\end{equation}
where $\mathcal{G}$ is a subclass of $F^{\varrho,r,q}_L$. Finally, the $n$th minimal error is defined as
\begin{equation}
	e_n^{(r)}(\mathcal{G},\mathcal{W},\delta_1,\delta_2)=\inf\limits_{\mathcal{A}\in\Phi_n}e^{(r)}(\mathcal{A},\mathcal{G},\mathcal{W},\delta_1,\delta_2).
\end{equation}
The aim is to develop an optimal algorithm and its efficient implementation by using GPUs.

Unless otherwise stated, all constants appearing in this paper (including those in the '$\mathcal{O}$', '$\Omega$', and '$\Theta$' notation) will only
depend on the parameters of the respective classes. Furthermore, the same symbol may be used for different constants.
\section{The Riemann-Maruyama quadrature for noisy information}\label{sec:alg}
Let $n\in\mathbb{N}$ and
\begin{equation}
	0=t_0<t_1<\ldots<t_n=T,
\end{equation}
be an arbitrary discretization on $[0,T]$. We denote by $\Delta t_i=t_{i+1}-t_i$, $i=0,1,\ldots,n-1$. We define the Riemann-Maruyama quadrature that use noisy evaluations of $X$ and $W$ by
\begin{equation}
	\mathcal{A}^{RM}_n(\tilde X,\tilde W)=\sum_{i=0}^{n-1}\tilde X(t_i)\cdot (\tilde W(t_{i+1})-\tilde W(t_{i})),
\end{equation}
where $(\tilde X,\tilde W)\in V_X(\delta_1)\times \mathcal{W}(\delta_2)$ for $\mathcal{W}\in\{\mathcal{W}_{s},\mathcal{\bar W}_{s},\mathcal{W}_{\alpha,\beta}\}$.
It is easy to see that the information cost of computing $\mathcal{A}^{RM}_n(\tilde X,\tilde W)$ is $2n$ noisy
evaluations of $X$ and $W$. The combinatory cost consists of $\mathcal{O}(n)$ arithmetic operations.

The aim of this section it to prove the following result.
\begin{theorem} 
\label{main_thm} Let us assume that $\varrho\in (0,1]$ and $r\geq 2$.
\begin{itemize}
\item [(i)] Let $s\geq 0$ and $q\in (r,+\infty)$. There exists a positive constant $C$, depending only on the parameters of the class $F^{\varrho,r,q}_L$ and $s$, such that for all $n\in\mathbb{N}$, $\delta_1,\delta_2\geq 0$, $X\in F^{\varrho,r,q}_L$, $(\tilde X, \tilde W)\in V_X(\delta_1)\times \mathcal{W}_{s}(\delta_2)$ it holds
\begin{equation}
\label{est_RMQ}
	\|\mathcal{I}(X,W)-\mathcal{A}^{RM}_n(\tilde X,\tilde W)\|_r\leq C\left(\max\limits_{0\leq i\leq n-1}(\Delta t_i)^{\varrho}+\delta_1+\delta_2+\delta_1\cdot\delta_2\right).
\end{equation}
\item [(ii)]Let $s\geq 0$ and $q\in (r,+\infty)$. There exists a positive constant $C$, depending only on the parameters of the class $F^{\varrho,r,q}_L$ and $s$, such that for all $n\in\mathbb{N}$, $\delta_1,\delta_2\geq 0$, $X\in F^{\varrho,r,q}_L$, $(\tilde X, \tilde W)\in V_X(\delta_1)\times \mathcal{\bar W}_{s}(\delta_2)$ it holds
\begin{eqnarray}
\label{est_RMQ3}
	\|\mathcal{I}(X,W)-\mathcal{A}^{RM}_n(\tilde X,\tilde W)\|_r&\leq& C\Bigl(\max\limits_{0\leq i\leq n-1}(\Delta t_i)^{\varrho}+\delta_1+\notag\\
&&\delta_2(1+\delta_1)\cdot(1+\sum\limits_{i=0}^{n-1}(\Delta t_i)^{1/2})\Bigr)
\end{eqnarray}
\item [(iii)] Let $\alpha,\beta\in (0,1]$ and $q\in [r,+\infty)$. There exists a positive constant $C$, depending only on the parameters of the class $F^{\varrho,r,q}_L$ and $\alpha$, $\beta$, such that for all $n\in\mathbb{N}$, $\delta_1,\delta_2\geq 0$, $X\in F^{\varrho,r,q}_L$, $(\tilde X, \tilde W)\in V_X(\delta_1)\times \mathcal{W}_{\alpha,\beta}(\delta_2)$ it holds
\begin{eqnarray}
\label{est_RMQ2}
	\|\mathcal{I}(X,W)-\mathcal{A}^{RM}_n(\tilde X,\tilde W)\|_r &\leq& C\Bigl(\max\limits_{0\leq i\leq n-1}(\Delta t_i)^{\varrho}+\delta_1+\notag\\
&&\delta_2(1+\delta_1)\cdot\sum\limits_{i=0}^{n-1}((\Delta t_i)^{\alpha}+(\Delta t_i)^{\beta/2})\Bigr).
\end{eqnarray}  
\end{itemize}  
\end{theorem}
\bf Proof. \rm Let $\tilde X\in V_X(\delta_1)$. We first show (\ref{est_RMQ}), where $\tilde W\in\mathcal{W}_s(\delta_2)$. Let the process $Z=\{Z(t)\}_{t\in [0,T]}$ be defined as $Z(t)=p_W(t,W(t))$. Then, by the It\^o formula we get that
\begin{equation}
	Z(t)=M(t)+V(t), \quad t\in [0,T],
\end{equation}
where

\begin{equation}
	V(t)=\int\limits_0^t\mathcal{L}p_W(z,W(z))\rd z,
\end{equation}
\begin{equation}
	M(t)=\int\limits_0^t\frac{\partial p_W}{\partial y}(z,W(z))\rd W(z).
\end{equation}
We stress that $\{V(t)\}_{t\in [0,T]}$ is continuous process with bounded variation, while $\{M(t)\}_{t\in [0,T]}$ is continuous martingale with respect to the filtration $\{\Sigma_t\}_{t\geq 0}$. Hence, the process $Z$ is continuous semimartingale. 

We denote by 
$$ \Delta Y_i = Y(t_{i+1})-Y(t_i) \qquad  i=0,1\ldots,n-1,
$$
for $Y\in\{W,Z\}$ and for all $t\in [0,T]$
\begin{equation}
	\hat X_n(t)=\sum_{i=0}^{n-1}X(t_i)\cdot\mathbf{1}_{(t_i,t_{i+1}]}(t),
\end{equation}
\begin{equation}
	\hat p_{X,n}(t)=\sum_{i=0}^{n-1}p_X(t_i,X(t_i))\cdot\mathbf{1}_{(t_i,t_{i+1}]}(t).
\end{equation}
Note that $\{X_n(t)\}_{t\in [0,T]}$ and $\{p_{X,n}(t)\}_{t\in [0,T]}$ are $\{\Sigma_t\}_{t\geq 0}$-progressively measurable simple processes.
Since $Z$ and $W$ are continuous semimartingales, by Property (v) at page 110 in \cite{legall} we can write the algorithm $\mathcal{A}^{RM}_n$ as follows
\begin{eqnarray}
\label{alg_main_decomp}
	\mathcal{A}^{RM}_n(\tilde X,\tilde W)&=&\sum\limits_{i=0}^{n-1}(X(t_i)+\delta_1\cdot p_X(t_i,X(t_i)))\cdot (\Delta W_i+\delta_2\cdot\Delta Z_i)\\
&=&\int\limits_0^T \hat X_n(t)\rd W(t)+\delta_1\int\limits_0^T \hat p_{X,n}(t)\rd W(t)\notag\\
\label{alg_decomp}
	&&+\delta_2\int\limits_0^T\hat X_n(t)\rd Z(t)+\delta_1\cdot\delta_2\int\limits_0^T \hat p_{X,n}(t)\rd Z(t).
\end{eqnarray}
We thus obtain
\begin{equation}
\label{err_est}
	\|\mathcal{I}(X,W)-\mathcal{A}^{RM}_n(\tilde X,\tilde W)\|_r\leq \sum_{i=1}^4 A_{i,n},
\end{equation}
where 
\begin{equation}
	A_{1,n}=\Biggl\|\int\limits_0^T \Bigl(\hat X_n(t)-X(t)\Bigr)\rd W(t)\Biggl\|_r,
\end{equation}
\begin{equation}
	A_{2,n}=\delta_1\cdot\Biggl\|\int\limits_0^T \hat p_{X,n}(t)\rd W(t)\Biggl\|_r,
\end{equation}
\begin{equation}
	A_{3,n}=\delta_2\cdot\Biggl\|\int\limits_0^T \hat X_n(t)\rd Z(t)\Biggl\|_r,
\end{equation}
\begin{equation}
	A_{4,n}=\delta_1\cdot\delta_2\cdot\Biggl\|\int\limits_0^T \hat p_{X,n}(t)\rd Z(t)\Biggl\|_r.
\end{equation}
By the Burkholder inequality and the H\"older continuity of $X$ in $r$th mean  we get
\begin{equation}
\label{a1_est}
	A_{1,n}\leq C_1\cdot\Biggl(\sum\limits_{i=0}^{n-1}\int\limits_{t_i}^{t_{i+1}}\mathbb{E}|X(t)-X(t_i)|^r \rd t\Biggr)^{1/r}\leq C_2\max\limits_{0\leq i\leq n-1}(\Delta t_i)^{\varrho},
\end{equation}
and
\begin{eqnarray}
\label{a2_est}
	A_{2,n}&\leq&\delta_1\cdot C_3\cdot\Biggl(\mathbb{E}\int\limits_0^T|\hat p_{X,n}(t)|^r\rd t\Biggr)^{1/r}=\delta_1\cdot C_3\cdot\Biggl(\sum_{i=0}^{n-1}\mathbb{E}|p_X(t_i,X(t_i))|^r\cdot\Delta t_i\Biggr)^{1/r}\notag\\
	&\leq& \delta_1\cdot C_4\Bigl(1+\Bigl\|\sup\limits_{t\in [0,T]}|X(t)|\Bigl\|_r\Bigr)\leq \delta_1\cdot C_4\Bigl(1+\Bigl\|\sup\limits_{t\in [0,T]}|X(t)|\Bigl\|_q\Bigr)\notag\\
&\leq& \delta_1\cdot C_5(1+L),
\end{eqnarray}
since $p_X$ is of at most linear growth. 
\newline
Since $Z=M+V$, from Definition 5.7 at page 109 in \cite{legall} we obtain
\begin{equation}
\label{a3_est_bb}
	A_{3,n}\leq \delta_2\cdot (B_{1,n}+B_{2,n}),
\end{equation}
where
\begin{equation}
	B_{1,n}=\Biggl\|\int\limits_0^T\hat X_n(t) \rd M(t)\Biggl\|_r,
\end{equation}
\begin{equation}
	B_{2,n}=\Biggl\|\int\limits_0^T\hat X_n(t) \rd V(t)\Biggl\|_r.
\end{equation}
Note that 
\begin{equation}
	\sup\limits_{t\in [0,T]}|\hat X_n(t)|\leq \sup\limits_{t\in [0,T]}|X(t)|,
\end{equation} 
\begin{equation}
	\sup\limits_{t\in [0,T]}|\hat p_{X,n}(t)|\leq 1+\sup\limits_{t\in [0,T]}|X(t)|,
\end{equation}
and
\begin{equation}
	\mathbb{E}\Bigl(\int\limits_0^T\Bigl|\frac{\partial p_W}{\partial y}(t,W(t))\Bigl|^r\rd t\Bigr)^{\frac{q}{q-r}}\leq C_6,
\end{equation}
\begin{equation}
	\mathbb{E}\Bigl(\int\limits_0^T|\mathcal{L}p_W(t,W(t))|^r\rd t\Bigr)^{\frac{q}{q-r}}\leq C_7,
\end{equation}
since $\frac{\partial p_W}{\partial y}$ and $\mathcal{L}p_W(t,y)$ are of at most linear growth. (The constants $C_6,C_7$  depend only on $T$, $s$, $r$, and $q$.) Hence, by the associativity property (see, for example, Property (ii) at page 109 in \cite{legall}), Burkholder and H\"older inequalities we get
\begin{eqnarray}
\label{b1_est}
	B_{1,n}&=&\Biggl\|\int\limits_0^T\hat X_n(t)\cdot\frac{\partial p_W}{\partial y}(t,W(t)) \rd W(t)\Biggl\|_r\leq C_8\Biggl(\mathbb{E}\int\limits_0^T|\hat X_n(t)|^r\cdot\Bigl|\frac{\partial p_W}{\partial y}(t,W(t))\Bigl|^r\rd t\Biggr)^{1/r}\notag\\
	&\leq& C_8\Bigl\|\sup\limits_{t\in [0,T]}|\hat X_n(t)|\Bigl\|_q\cdot\Biggl(\mathbb{E}\Bigl(\int\limits_0^T\Bigl|\frac{\partial p_W}{\partial y}(t,W(t))\Bigl|^r\rd t\Bigr)^{\frac{q}{q-r}}\Biggr)^{\frac{1}{r}-\frac{1}{q}}\notag\\
&\leq& C_9\Bigl\|\sup\limits_{t\in [0,T]}|X(t)|\Bigl\|_q\leq C_9L,
\end{eqnarray}
and
\begin{eqnarray}
\label{b2_est}
	B_{2,n}&=&\Biggl\|\int\limits_0^T\hat X_n(t)\cdot\mathcal{L}p_W(t,W(t)) \rd t\Biggl\|_r\leq C_{10}\Biggl(\mathbb{E}\int\limits_0^T|\hat X_n(t)|^r\cdot|\mathcal{L}p_W(t,W(t))|^r\rd t\Biggr)^{1/r}\notag\\
&\leq& C_{11} \Bigl\|\sup\limits_{t\in [0,T]}|\hat X_n(t)|\Bigl\|_q\cdot\Biggl(\mathbb{E}\Bigl(\int\limits_0^T|\mathcal{L}p_W(t,W(t))|^r\rd t\Bigr)^{\frac{q}{q-r}}\Biggr)^{\frac{1}{r}-\frac{1}{q}}\notag\\
&\leq& C_{12}\Bigl\|\sup\limits_{t\in [0,T]}|X(t)|\Bigl\|_q\leq C_{12}L,
\end{eqnarray}
where $C_9,C_{12}$ depend only on $T$, $s$, $q$, $r$ and $L$. Therefore, by (\ref{b1_est}), (\ref{b2_est}), and (\ref{a3_est_bb}) we arrive at
\begin{equation}
\label{a3_est}
	A_{3,n}\leq C_{13}\cdot\delta_2.
\end{equation}
By proceeding analogously as for $A_{3,n}$ we obtain 
\begin{equation}
\label{a4_est}
	A_{4,n}\leq C_{14}\cdot\delta_1\cdot\delta_2.
\end{equation} 
Combining (\ref{err_est}), (\ref{a1_est}), (\ref{a2_est}), (\ref{a3_est}), and (\ref{a4_est}) we get (\ref{est_RMQ}), which ends the proof of (\ref{est_RMQ}).

We now justify (\ref{est_RMQ2}) and (\ref{est_RMQ3}). In this cases the process $Z$ is not necessarily a semimartingale. Hence, we use the following decomposition of $\mathcal{A}_n^{RM}$, that follows directly from (\ref{alg_main_decomp}),
\begin{eqnarray}
	&&\mathcal{A}^{RM}_n(\tilde X,\tilde W)=\int\limits_0^T \hat X_n(t)\rd W(t)+\delta_1\int\limits_0^T \hat p_{X,n}(t)\rd W(t)\notag\\
	&&+\delta_2\sum\limits_{i=0}^{n-1}X(t_i)\cdot \Delta Z_i+\delta_1\cdot\delta_2\sum\limits_{i=0}^{n-1}p_X(t_i,X(t_i))\cdot\Delta Z_i.
\end{eqnarray} 
We have that
\begin{equation}
\label{err_est2}
	\|\mathcal{I}(X,W)-\mathcal{A}^{RM}_n(\tilde X,\tilde W)\|_2\leq \sum_{i=1}^4 D_{i,n},
\end{equation}
where 
\begin{equation}
\label{est_d1}
	D_{1,n}=A_{1,n},
\end{equation}
\begin{equation}
\label{est_d2}
	D_{2,n}=A_{2,n},
\end{equation}
\begin{equation}
\label{est_d3}
	D_{3,n}=\delta_2\cdot\Biggl\|\sum\limits_{i=0}^{n-1}X(t_i)\cdot\Delta Z_i\Biggl\|_r,
\end{equation}
\begin{equation}
\label{est_d4}
	D_{4,n}=\delta_1\cdot\delta_2\cdot\Biggl\|\sum\limits_{i=0}^{n-1}p_X(t_i,X(t_i))\cdot\Delta Z_i\Biggl\|_r.
\end{equation}
For $D_{1,n}$ and $D_{2,n}$ we use the bounds (\ref{a1_est}), (\ref{a2_est}) obtained for $A_{1,n}$ and $A_{2,n}$, respectively. However, for $D_{3,n}$ and $D_{4,n}$ we have to differ between the case when $\tilde W\in\mathcal{\bar W}_s(\delta_2)$ and $\tilde W\in\mathcal{W}_{\alpha,\beta}(\delta_2)$.

Let $\tilde W\in\mathcal{\bar W}_s(\delta_2)$. Since in this case $p_W\in\mathcal{\bar K}^2_s$ we get, by the mean value theorem, that 

\begin{equation}
	\label{local_Lip_pw}
		|p_W(t,x)-p_W(z,y)|\leq C_1\Bigl(|t-z|+(1+|x|^s+|y|^s)\cdot |x-y|\Bigr), 
\end{equation}
for all $t,z\in [0,T]$ and $x,y\in\mathbb{R}$, where $\bar C>0$ depends only on $T$ and $s$. This implies
\begin{equation}
	|\Delta Z_i|^{r\gamma}\leq C_2(\Delta t_i)^{r\gamma}+C_3(1+|W(t_i)|^{rs\gamma}+|W(t_{i+1})|^{rs\gamma})\cdot|\Delta W_i|^{r\gamma},
\end{equation}
for $i=0,1,\ldots,n-1$, where $\gamma=q/(q-r)$. Hence, by the H\"older inequality
\begin{equation}
	\| X(t_i)\cdot \Delta Z_i\|_r\leq \Bigl\|\sup\limits_{t\in [0,T]}|X(t)| \Bigl\|_{q}\cdot\Bigl(\mathbb{E}|\Delta Z_i|^{r\gamma}\Bigr)^{1/(r\gamma)},
\end{equation}
and
\begin{eqnarray}
	\mathbb{E}|\Delta Z_i|^{r\gamma}&\leq& C_3\Bigl(\mathbb{E}(1+|W(t_i)|^{rs\gamma}+|W(t_{i+1})|^{rs\gamma})^2\Bigr)^{1/2}\cdot\Bigl(\mathbb{E}|\Delta W_i|^{2r\gamma}\Bigr)^{1/2}\notag\\
&&+C_2(\Delta t_i)^{r\gamma}\leq C_2(\Delta t_i)^{r\gamma}+C_4(\Delta t_i)^{r\gamma/2},
\end{eqnarray}
since $W$ has all absolute moments bounded. Therefore,
\begin{equation}
	\| X(t_i)\cdot \Delta Z_i\|_r\leq C_5 \Bigl\|\sup\limits_{t\in [0,T]}|X(t)| \Bigl\|_{q}\cdot\Bigl(\Delta t_i+(\Delta t_i)^{1/2}\Bigr),
\end{equation}
which implies that
\begin{equation}
D_{3,n}=\delta_2\cdot\Biggl\|\sum\limits_{i=0}^{n-1}X(t_i)\cdot \Delta Z_i\Biggl\|_r\leq\delta_2\cdot\sum\limits_{i=0}^{n-1}\|X(t_i)\cdot \Delta Z_i\|_r\leq C_6\cdot\delta_2\cdot(T+\sum\limits_{i=0}^{n-1}(\Delta t_i)^{1/2}).
\end{equation}
For $D_{4,n}$ we proceed analogously as for $D_{3,n}$.  This gives (\ref{est_RMQ3}).

Finally, let $\tilde W \in\mathcal{W}_{\alpha,\beta}(\delta_2)$. Since $X(t_i)$ and $\Delta W_i$ are independent, we have in this case that
\begin{eqnarray}
	D_{3,n}=\delta_2\cdot\Biggl\|\sum\limits_{i=0}^{n-1}X(t_i)\cdot \Delta Z_i\Biggl\|_r
	&\leq&\delta_2\cdot\sum\limits_{i=0}^{n-1} \| |X(t_i)| \|_r\cdot\|(\Delta t_i)^{\alpha}+|\Delta W_i|^{\beta}\|_r\notag\\
&\leq& \delta_2\cdot\Bigl\|\sup\limits_{t\in [0,T]}|X(t)|\Bigl\|_r\cdot\sum\limits_{i=0}^{n-1} \left((\Delta t_i)^{\alpha}+\||\Delta W_i|^{\beta}\|_r\right)\notag\\
&\leq& \delta_2\cdot L\cdot\sum\limits_{i=0}^{n-1} \left((\Delta t_i)^{\alpha}+m^{\beta}_{r\beta}(\Delta t_i)^{\beta/2}\right),
\end{eqnarray}
where $m_{r\beta}=\|Z\|_{r\beta}$ and $Z$ is a standard normal random variable with mean zero and variance equal to $1$. 
For $D_{4,n}$ we proceed analogously as for $D_{3,n}$. This ends the proof.  \ \ \ $\blacksquare$
\newline\newline
Directly from Theorem \ref{main_thm} we have the following corollary that states the worst-case error of the algorithm $\mathcal{A}^{RM}_n$  in the class $F^{\varrho,r,q}_L$.
\begin{corollary}
\label{cor_RM_upp} Let $\varrho\in (0,1]$, $r\geq 2$, and let us consider the Riemann-Maruyama quadrature $\mathcal{A}^{RM}_n$ based on the equidistant mesh $t_i=iT/n$, $i=0,1,\ldots,n$. 
\begin{itemize}
	\item [(i)] Let $s\geq 0$ and $q\in (r,+\infty)$. Then
				\begin{equation}
					\label{upp_RM_class1}
					e^{(r)}(\mathcal{A}^{RM}_n,F^{\varrho,r,q}_L,\mathcal{W}_{s},\delta_1,\delta_2)=\mathcal{O}(n^{-\varrho}+\delta_1+\delta_2),
				\end{equation}
				as $n\to +\infty$ and $\max\{\delta_1,\delta_2\}\to 0+$.
	\item [(ii)] Let $s\geq 0$ and $q\in (r,+\infty)$. Then
				\begin{equation}
					\label{upp_RM_class2}
					e^{(r)}(\mathcal{A}^{RM}_n,F^{\varrho,r,q}_L,\mathcal{\bar W}_{s},\delta_1,\delta_2)=\mathcal{O}(n^{-\varrho}+\delta_1+\delta_2\cdot (1+n^{1/2})),
				\end{equation}
				as $n\to +\infty$ and $\max\{\delta_1,\delta_2\}\to 0+$.			
	\item [(iii)] Let $\alpha,\beta\in (0,1]$ and $q\in [r,+\infty)$. Then
	\begin{equation}
					\label{upp_RM_class3}
					e^{(r)}(\mathcal{A}^{RM}_n,F^{\varrho,r,q}_L,\mathcal{W}_{\alpha,\beta},\delta_1,\delta_2)=\mathcal{O}(n^{-\varrho}+\delta_1+\delta_2\cdot n^{1-\min\{\alpha,\beta/2\}}),
				\end{equation}
	as $n\to +\infty$ and $\max\{\delta_1,\delta_2\}\to 0+$.
\end{itemize}	
\end{corollary}
Let us comment on the result obtained so far.
\begin{remark} As we can see from Theorem \ref{main_thm} and Corollary \ref{cor_RM_upp} domination of the noise term become more on more visible as the regularity of  disturbing functions $p_W$ is decreasing. 
\end{remark}
\begin{remark}
	\label{rem_1}
	We considered the setting that we called analytic noise model, since we assumed certain form  of the noise via disturbance function $p$. Of course another approach is possible. Namely, one can assume that the exact values of $X$ are corrupted by noise in the following way
	\begin{equation}
		\tilde X(t_i)=X(t_i)+\varepsilon_i,
	\end{equation}
	where  $(\varepsilon_i)_{i=0,1,\ldots,n}$ are $\sigma\Bigl(\bigcup_{t\geq 0}\Sigma_t\Bigr)$-measurable random variables. Preliminary estimates indicate that it is possible to achieve upper bounds like in Theorem \ref{main_thm}, under certain assumptions on the discrete-time process $(\varepsilon_i)_{i=0,1,\ldots,n}$. We postpone this problem to our future work.
\end{remark}
\section{Lower bounds and optimality of the Riemann-Maruyama quadrature} \label{sec:low}
In this section we investigate lower bounds on the worst-case error of an arbitrary
algorithm from $\Phi_n$  and, in particular cases, we establish optimality of the Riemann-Maruyama algorithm $\mathcal{A}^{RM}_n$. We concentrate on the class $\mathcal{W}_s$ of noisy evaluations of $W$. Essentially sharp lower bounds in the classes $\mathcal{\bar W}_s$ and
 $\mathcal{W}_{\alpha,\beta}$ are left as an open problem.
 
The following lemma follows directly from (91) in \cite{MoPr17}, where the lower bound on the error for approximating It\^o integrals of deterministic functions from the H\"older class has been established.
\begin{lemma} 
\label{lower_bounds}
Let $\varrho\in (0,1]$, $r\geq 2$, $q\in (r,+\infty)$, and $s\geq 0$. Then 
\begin{equation}
	\label{low_b_inex}
	e^{(r)}(F^{\varrho,r,q}_L,\mathcal{W}_s,\delta_1,\delta_2) = \Omega\Bigl( \max\{n^{-\varrho}, \delta_1\}\Bigr).
\end{equation}
\end{lemma}
From Corollary \ref{cor_RM_upp} and Lemma \ref{lower_bounds} we get the  main result of this paper.
\begin{theorem} 
\label{OPT_RME_ST}
Let $\varrho\in (0,1]$, $r\geq 2$, $q\in (r,+\infty)$, and $s\geq 0$. Then  the $n$th minimal error satisfies
\begin{equation}
		\label{OPT_B_NTH_ME1}
		e^{(r)}_n(F^{\varrho,r,q}_L,\mathcal{W}_{s},\delta_1,0)=\Theta(n^{-\varrho}+\delta_1),
	\end{equation}
and
	\begin{equation}
		\label{OPT_B_NTH_ME2}
		e^{(r)}_n(F^{\varrho,r,q}_L,\mathcal{W}_{s},\delta_1,\delta_1)=\Theta(n^{-\varrho}+\delta_1),
	\end{equation}
	as $n\to +\infty$ and $\delta_1\to 0+$.
An optimal algorithm is the Riemann-Maruyama quadrature $\mathcal{A}^{RM}_n$ based on the equidistant discretization $t_i=iT/n$, $i=0,1,\ldots,n$. 
\end{theorem}
The results above hold for particular values of the precision parameter $\delta_2$, namely, for $\delta_2=0$ and $\delta_2=\delta_1$. In general case preliminary estimates suggest that in order to establish dependence of the lower bounds also on $\delta_2$ completely new technique is required. (We stress that the results from \cite{Pla96} are not applicable here, since we consider a different model of noise.) Nevertheless, for the algorithm $\mathcal{A}^{RM}_n$ we have the following sharp (worst-case) error bounds in the case of arbitrary $\delta_2$.
\begin{proposition} \label{prop_lower} 
Let $\varrho\in (0,1]$, $r\geq 2$, $q\in (r,+\infty)$, $s\geq 0$, and let us consider the Riemann-Maruyama quadrature $\mathcal{A}^{RM}_n$ based on the equidistant mesh $t_i=iT/n$, $i=0,1,\ldots,n$.
	 Then
				\begin{equation}
					\label{opt_RM_class}
					e^{(r)}(\mathcal{A}^{RM}_n,F^{\varrho,r,q}_L,\mathcal{W}_{s},\delta_1,\delta_2)=\Theta(n^{-\varrho}+\delta_1+\delta_2),
				\end{equation}
	as $n\to +\infty$ and $\max\{\delta_1,\delta_2\}\to 0+$.			
\end{proposition}
\noindent
{\bf Proof.} Upper bounds in (\ref{opt_RM_class}) follows directly from Corollary \ref{cor_RM_upp}. In the case when $\delta_1\geq 0=\delta_2$, the lower bound $\Omega(\max\{n^{-\varrho},\delta_1\})$ again follows from (91) in \cite{MoPr17}.

Now we consider the case $\delta_2\geq 0=\delta_1$. Let us take
\begin{equation}
	X_0(t,\omega)=L, \quad (t,\omega)\in [0,T]\times\Omega,
\end{equation}
and take $\tilde W\in\mathcal{W}_s(\delta_2)$ of the following form
\begin{equation}
	\tilde W(t)=W(t)+\delta_2\cdot t.
\end{equation}
We get that
\begin{equation}
	\mathcal{I}(X_0,W)=L\cdot W(T),
\end{equation}
and
\begin{equation}
	\mathcal{A}^{RM}_n(X_0,\tilde W)=L\sum\limits_{i=0}^{n-1}\Delta W_i+L\delta_2\sum\limits_{i=0}^{n-1}\Delta t_i=\mathcal{I}(X_0,W)+L\delta_2 T,
\end{equation}
which gives 
\begin{equation}
	\|\mathcal{I}(X_0,W)-\mathcal{A}^{RM}_n(X_0,\tilde W)\|_r=LT\delta_2.
\end{equation}
This implies the thesis. \ \ \ $\blacksquare$
\begin{remark}\rm In the case of exact information (i.e., $\delta_1=\delta_2=0$)  we know, by the results of \cite{hein1}, that even randomized adaptive information does not help, and the rate $n^{-\varrho}$ is optimal.
\end{remark}
\section{Numerical results} \label{sec:num}
We present results for the  Riemann-Maruyama quadrature $\mathcal{A}_n^{RM}$. There will be four exemplary problems presented, where for the first one we know the exact solution and for the others  we need to assume some convergence of the analyzed method to estimate the obtained error. In the end of this section some practical guides on how to implement the algorithm efficiently using GPUs will be presented together with the discussion about the obtainable speedup of using such architecture.

\begin{figure}[h]
\centering
\subfloat[][$p_X(t,x)=xt^2\; p_W(t,x)=xt^2$]{
\includegraphics[width=0.45\linewidth]{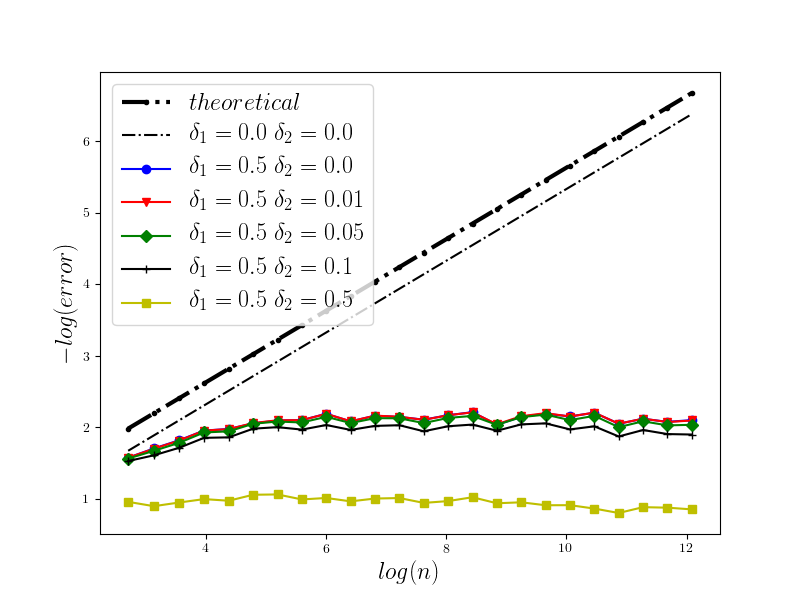}}
\quad
\subfloat[][$p_X(t,x)=xt^2\; p_W(t,x)=xt^2$]{
\includegraphics[width=0.45\linewidth]{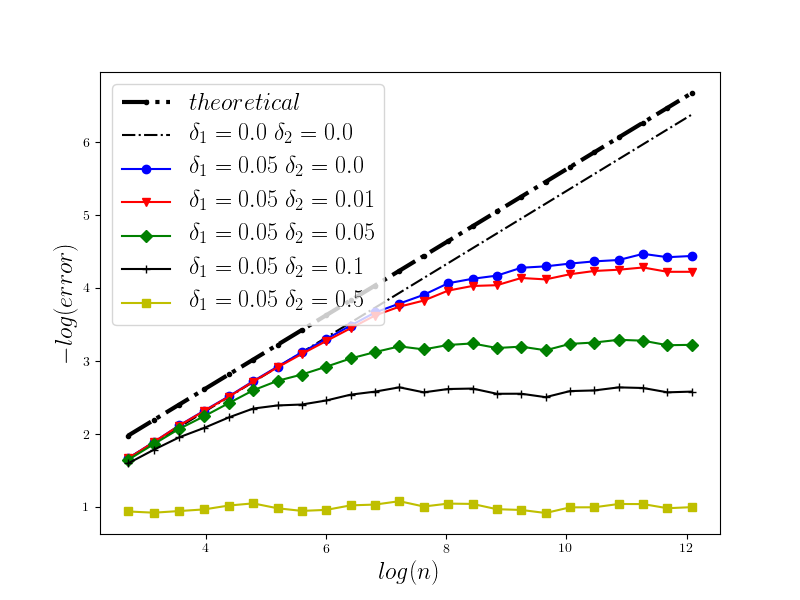}}
 \\
\subfloat[][$p_X(t,x)=x\; p_W(t,x)=1$]{
\includegraphics[width=0.45\linewidth]{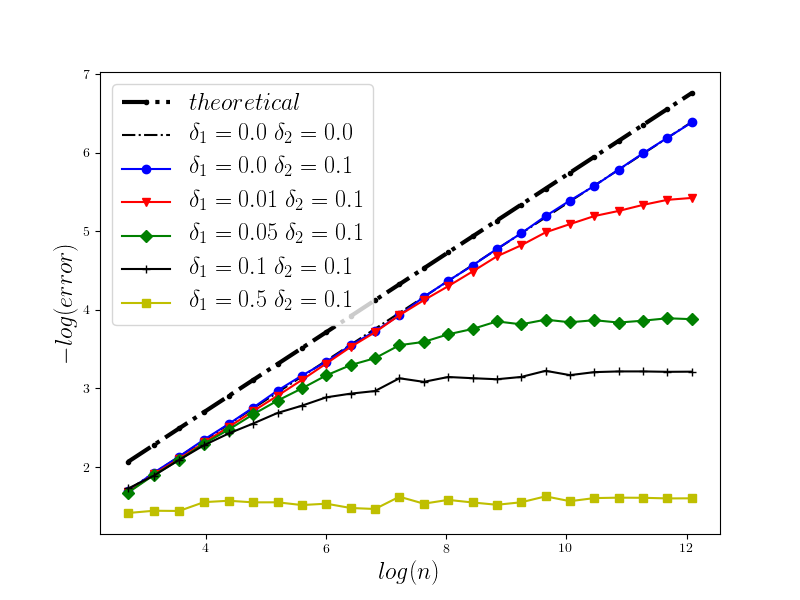}}
\quad
\subfloat[][$p_X(t,x)=1\; p_W(t,x)=x$]{
\includegraphics[width=0.45\linewidth]{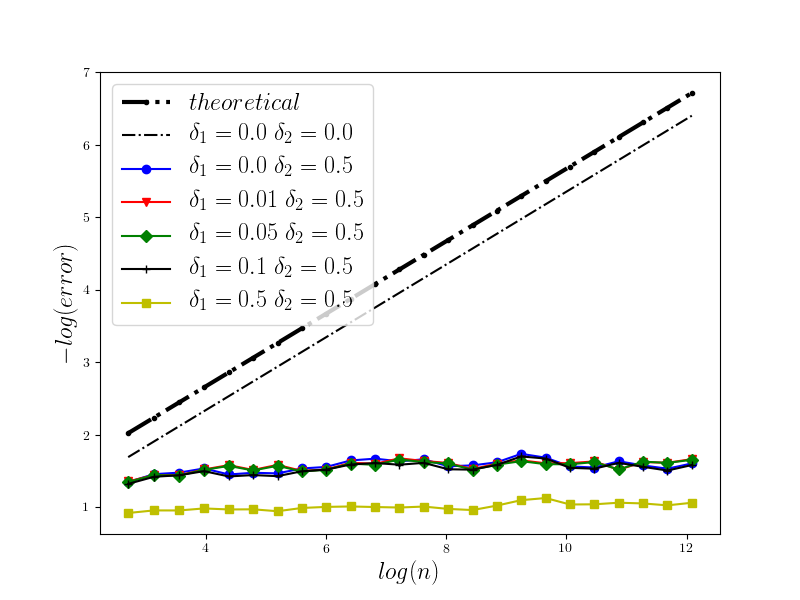}}

\caption{Error behavior for Riemann-Maruyama quadrature under exact/inexact information for problem \eqref{NR_PROB_1}.
}

\label{fig1}
\end{figure}


\begin{figure}
\centering
\subfloat[][$p_X(t,x)=1\; p_W(t,x)=x$]{
\includegraphics[width=0.45\linewidth]{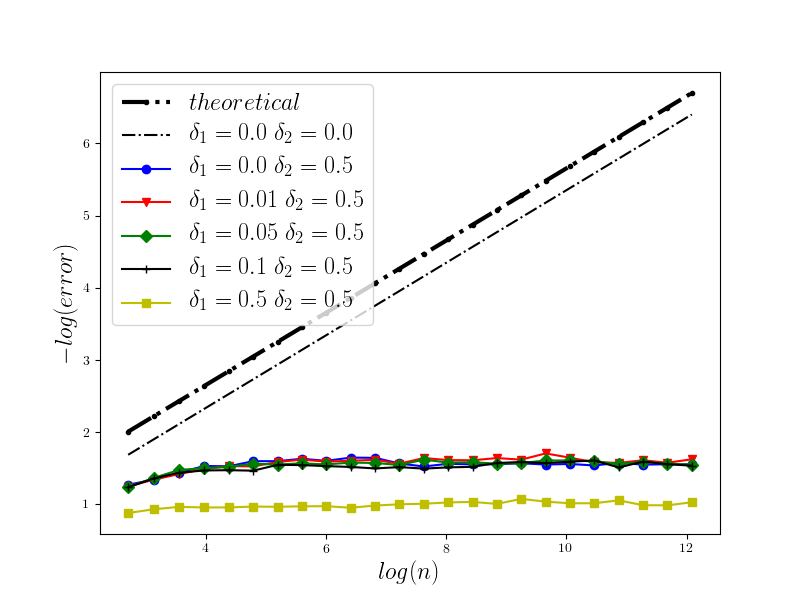}}
\quad
\subfloat[][$p_X(t,x)=1\; p_W(t,x)=x$]{
\includegraphics[width=0.45\linewidth]{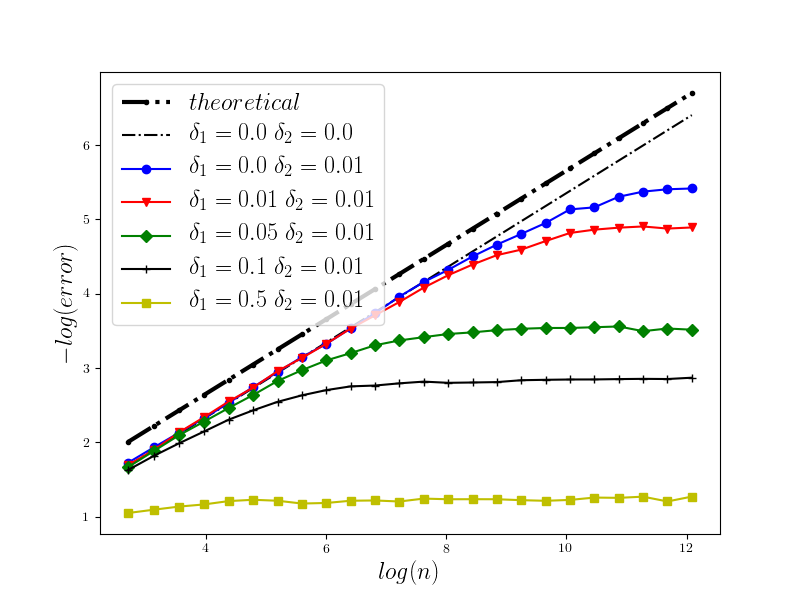}}
 \\
\subfloat[][$p_X(t,x)=1\; p_W(t,x)=x$]{
\includegraphics[width=0.45\linewidth]{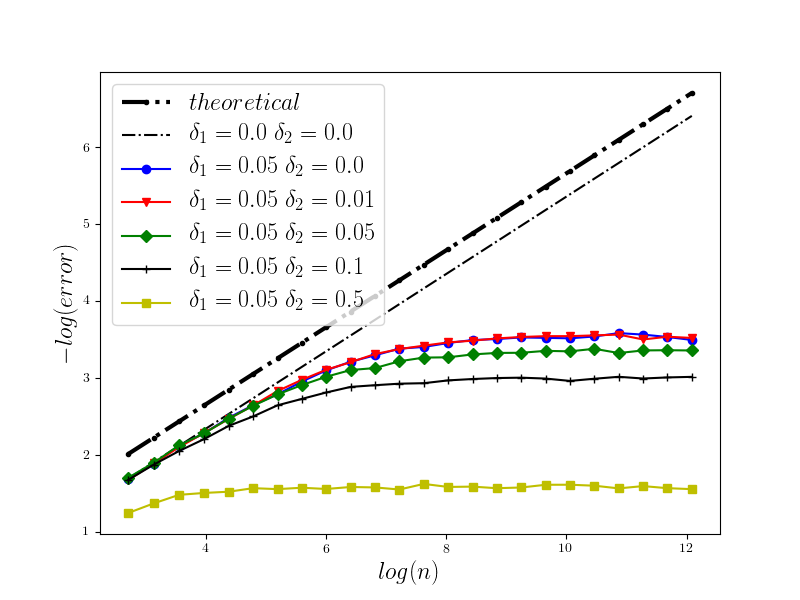}}
\quad
\subfloat[][$p_X(t,x)=1\; p_W(t,x)=x$]{
\includegraphics[width=0.45\linewidth]{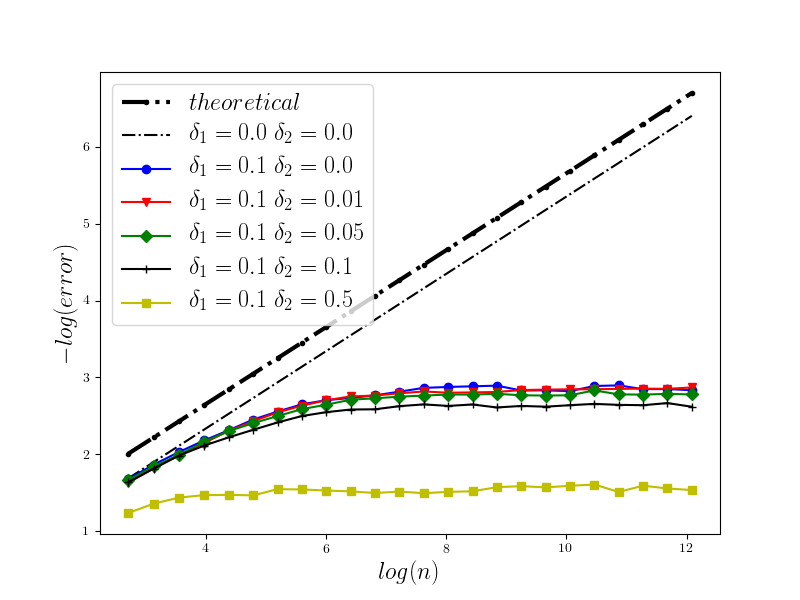}}

\caption{Error behavior for Riemann-Maruyama quadrature under exact/inexact information for problem \eqref{NR_PROB_2}.
}
\label{fig2}
\end{figure}

\begin{figure}
\centering
\subfloat[][$p_X(t,x)=1\; p_W(t,x)=x$]{
\includegraphics[width=0.45\linewidth]{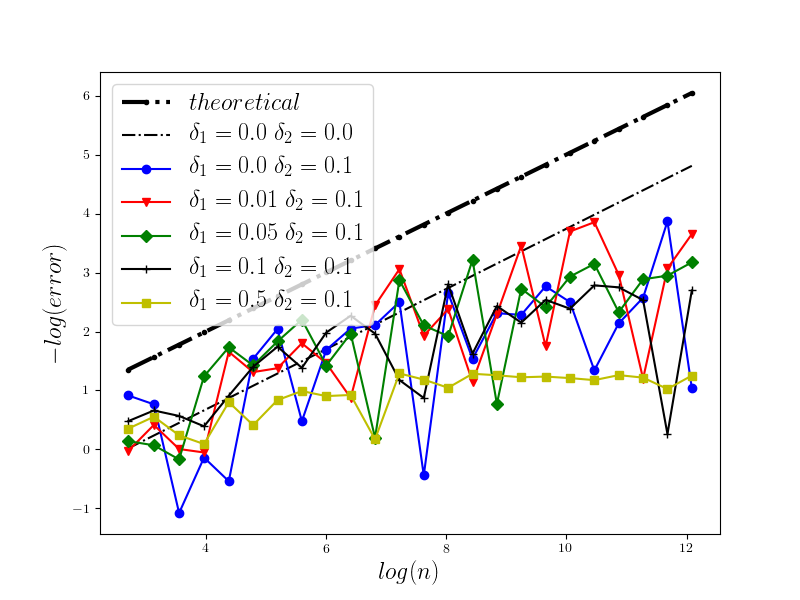}}
\quad
\subfloat[][$p_X(t,x)=1\; p_W(t,x)=x$]{
\includegraphics[width=0.45\linewidth]{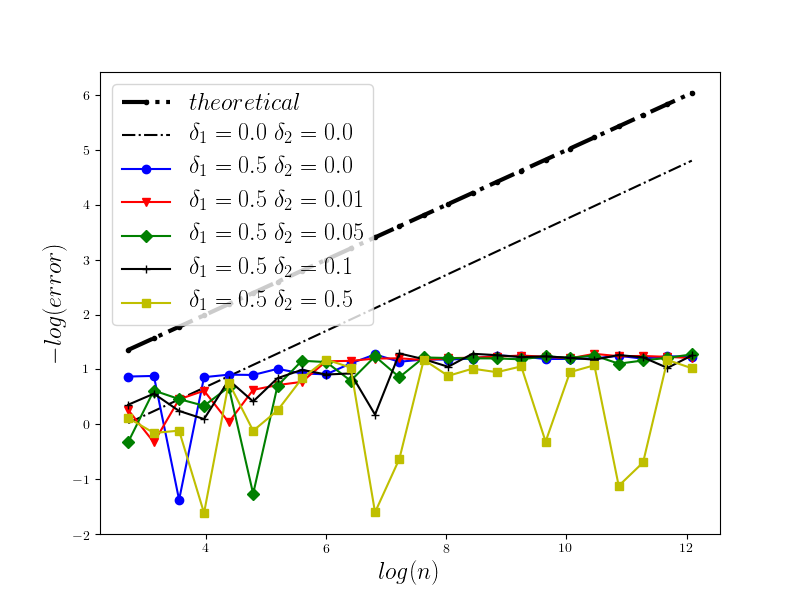}}
 \\
\subfloat[][$p_X(t,x)=1\; p_W(t,x)=x$]{
\includegraphics[width=0.45\linewidth]{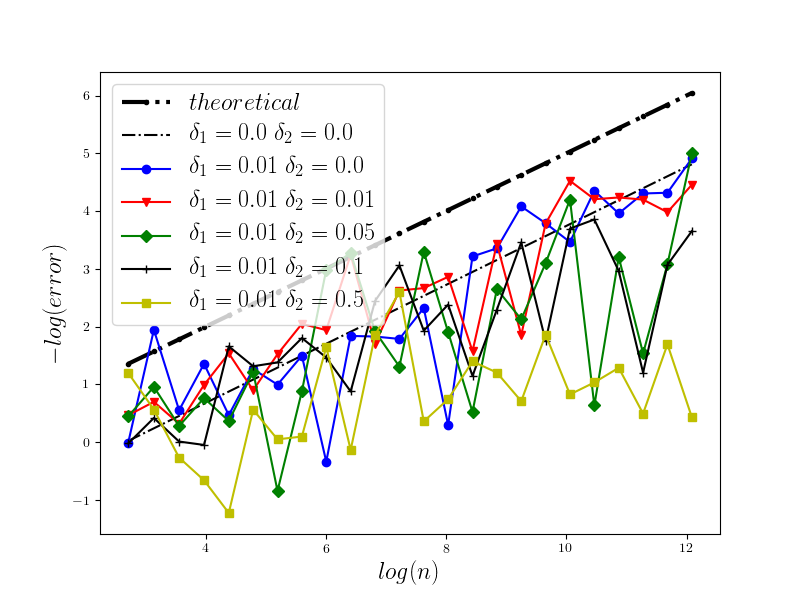}}
\quad
\subfloat[][$p_X(t,x)=1\; p_W(t,x)=x$]{
\includegraphics[width=0.45\linewidth]{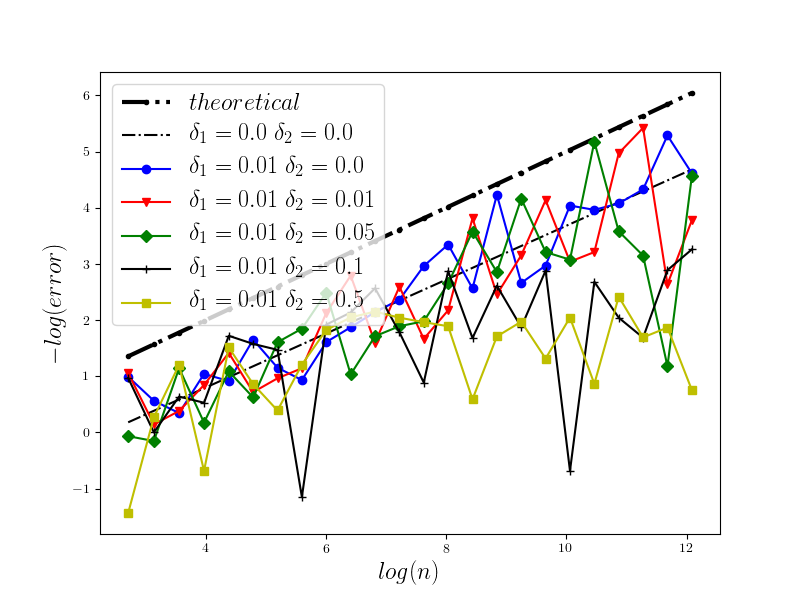}}

\caption{Error behavior for Riemann-Maruyama quadrature under exact/inexact information for problem \eqref{NR_PROB_3}.
}
\label{fig3}
\end{figure}

\begin{figure}
\centering
\subfloat[][$p_X(t,x)=xt^2\; p_W(t,x)=1$]{
\includegraphics[width=0.45\linewidth]{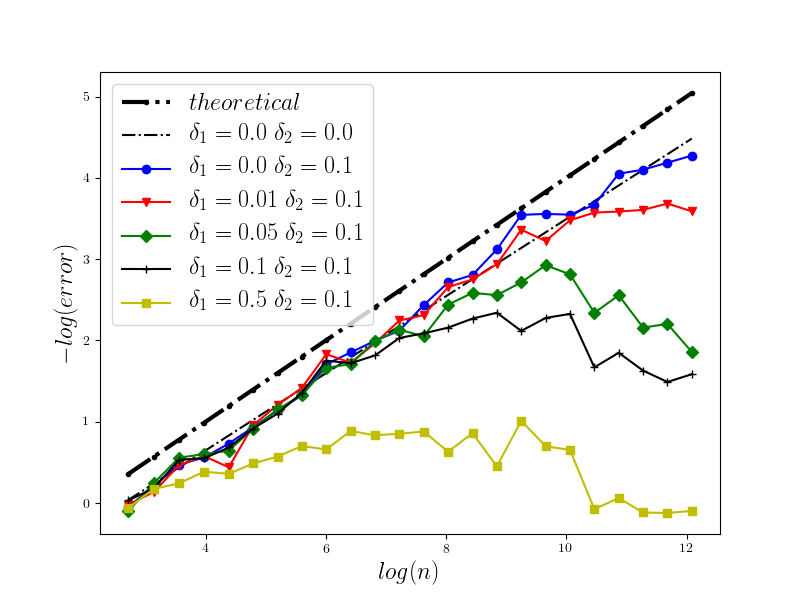}}
\quad
\subfloat[][$p_X(t,x)=xt^2\; p_W(t,x)=xt^2$]{
\includegraphics[width=0.45\linewidth]{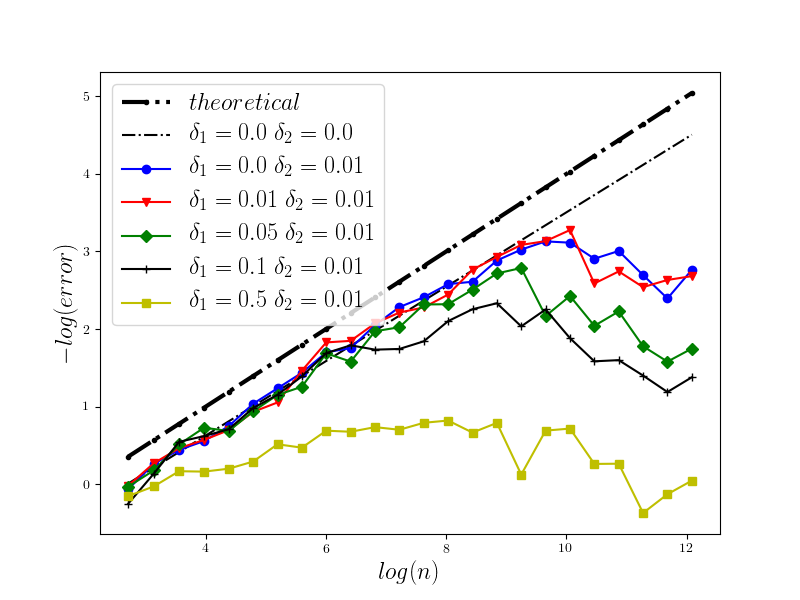}}
 \\
\subfloat[][$p_X(t,x)=1\; p_W(t,x)=x$]{
\includegraphics[width=0.45\linewidth]{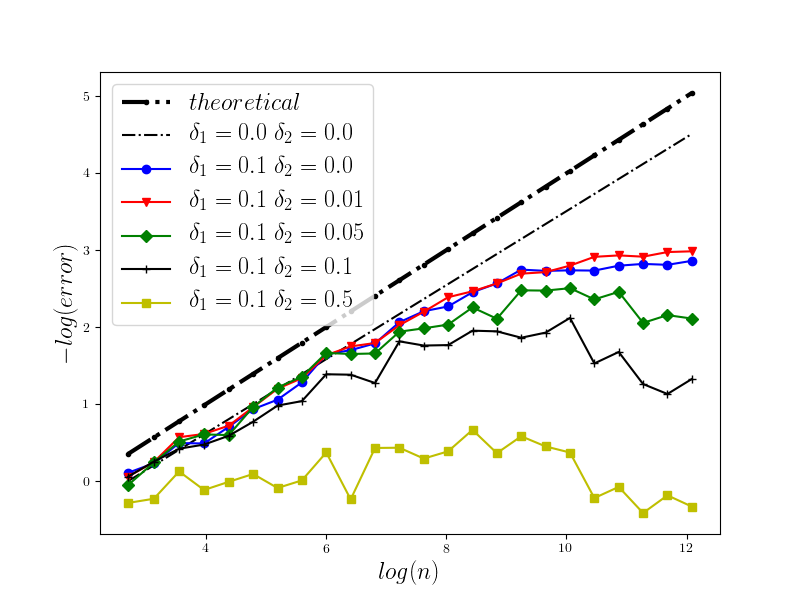}}
\quad
\subfloat[][$p_X(t,x)=xt^2\; p_W(t,x)=xt^2$]{
\includegraphics[width=0.45\linewidth]{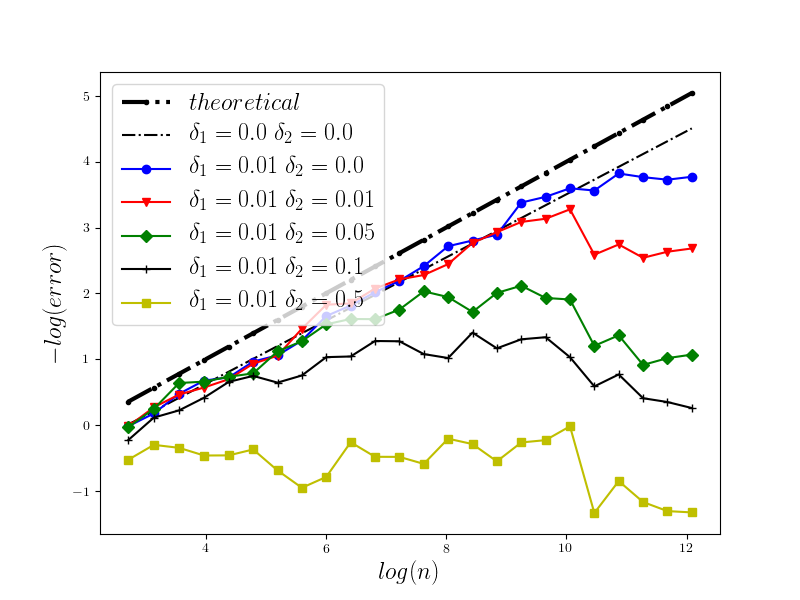}}

\caption{Error behavior for Riemann-Maruyama quadrature under exact/inexact information for problem \eqref{NR_PROB_4}.
}
\label{fig4}
\end{figure}

\begin{figure}
\centering
\subfloat[][$p_X(t,x)=1\; p_W(t,x)=x$]{
\includegraphics[width=0.45\linewidth]{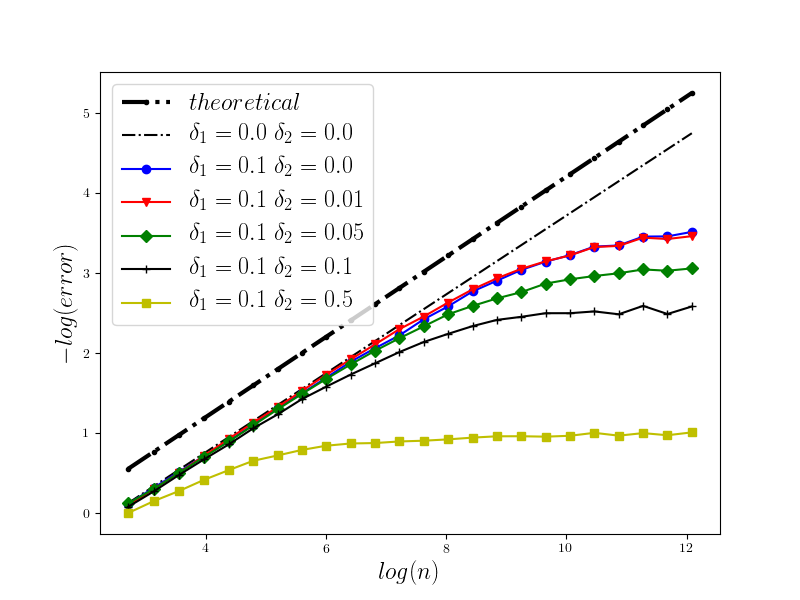}}
\quad
\subfloat[][$p_X(t,x)=xt^2\; p_W(t,x)=1$]{
\includegraphics[width=0.45\linewidth]{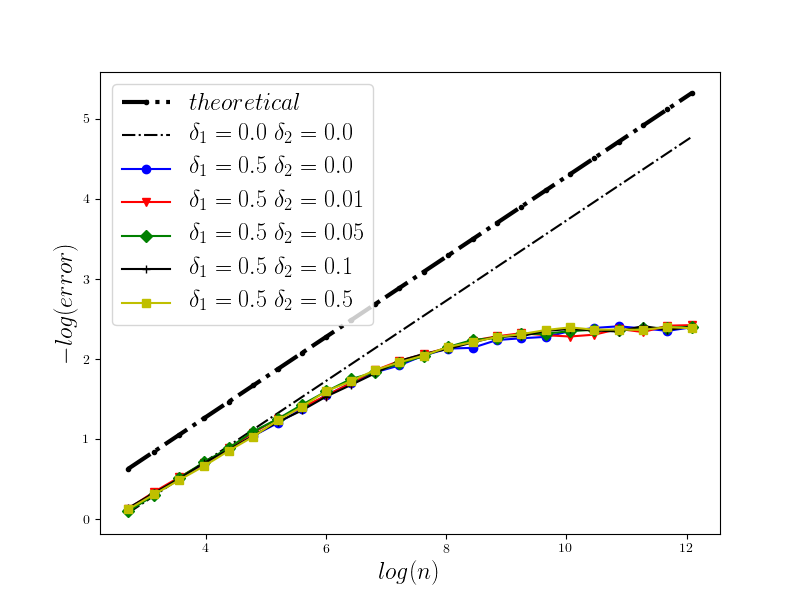}}
 \\
\subfloat[][$p_X(t,x)=xt^2\; p_W(t,x)=1$]{
\includegraphics[width=0.45\linewidth]{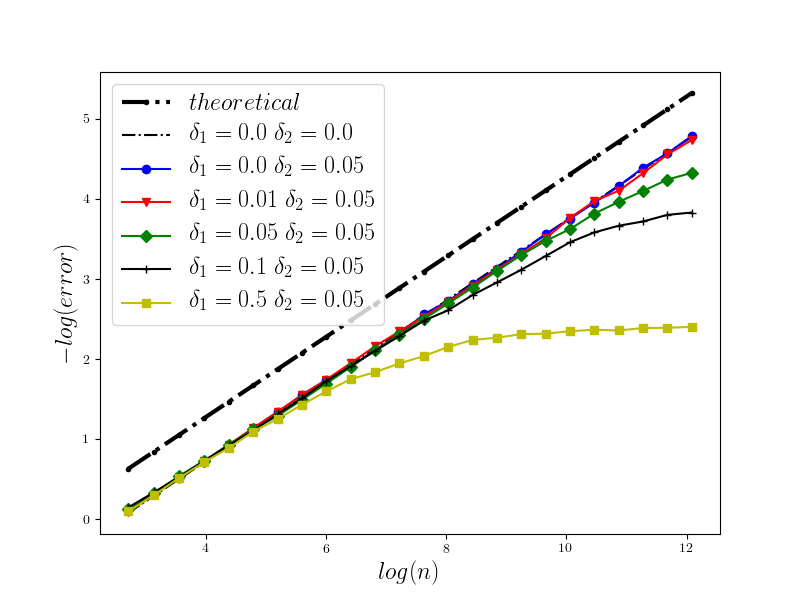}}
\quad
\subfloat[][$p_X(t,x)=x\; p_W(t,x)=1$]{
\includegraphics[width=0.45\linewidth]{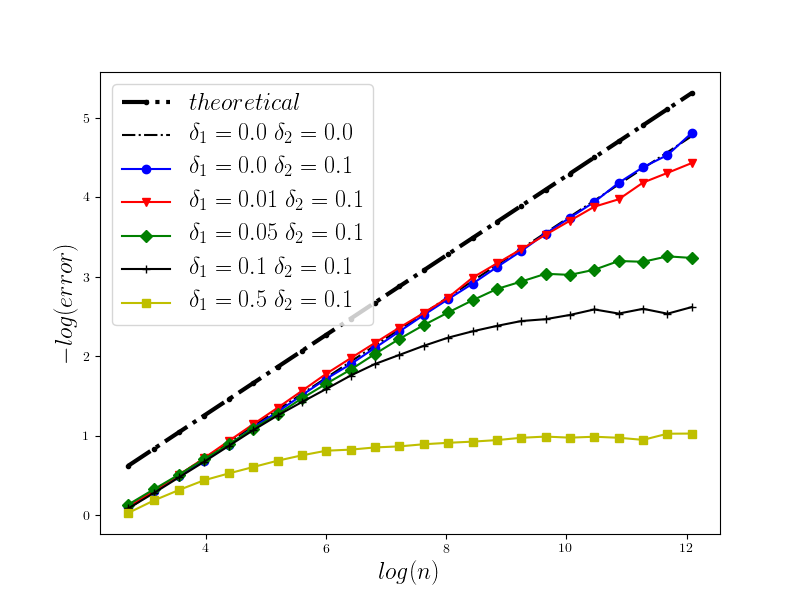}}

\caption{Behavior of weak error under exact/inexact information for problem \eqref{weak_app_Y}.
}
\label{fig6}
\end{figure}

\begin{figure}
\centering
\subfloat[][Problem \eqref{NR_PROB_1}]{
\includegraphics[width=0.45\linewidth]{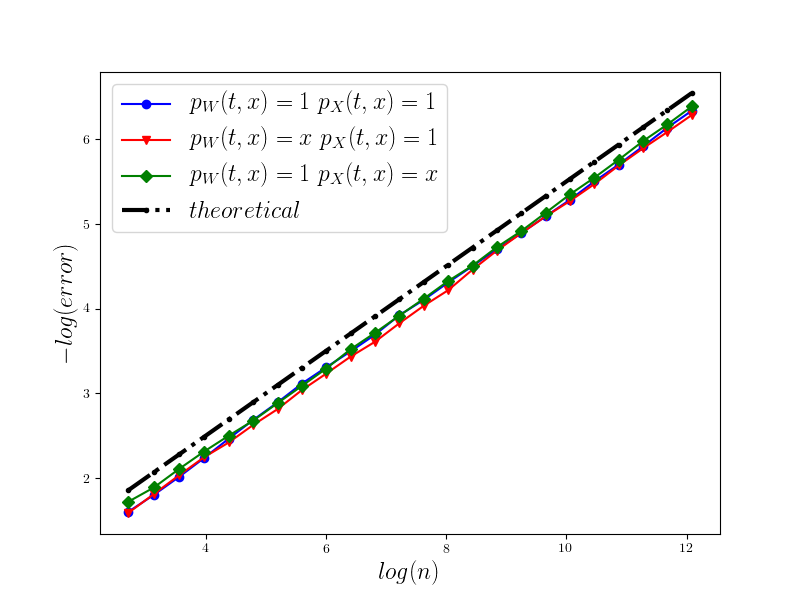}}
\quad
\subfloat[][Problem \eqref{NR_PROB_2}]{
\includegraphics[width=0.45\linewidth]{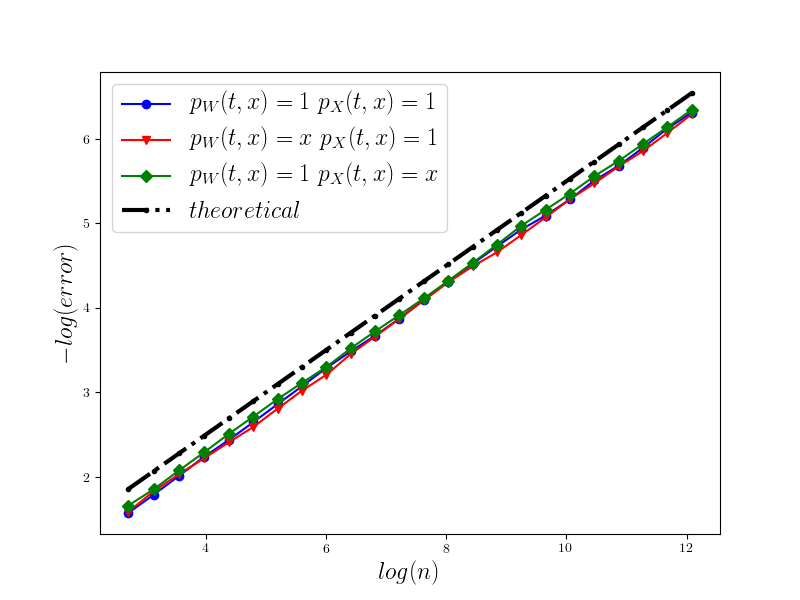}}
 \\
\subfloat[][Problem \eqref{NR_PROB_3}]{
\includegraphics[width=0.45\linewidth]{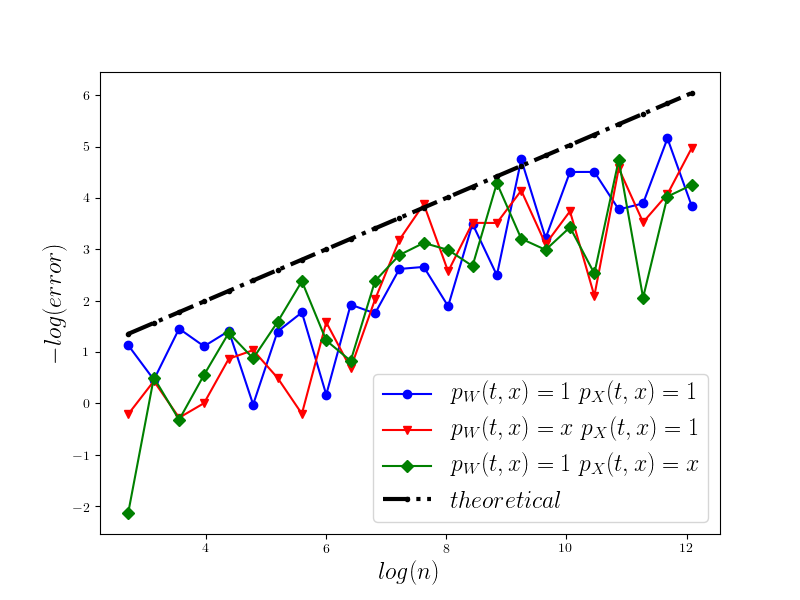}}
\quad
\subfloat[][Problem \eqref{NR_PROB_4}]{
\includegraphics[width=0.45\linewidth]{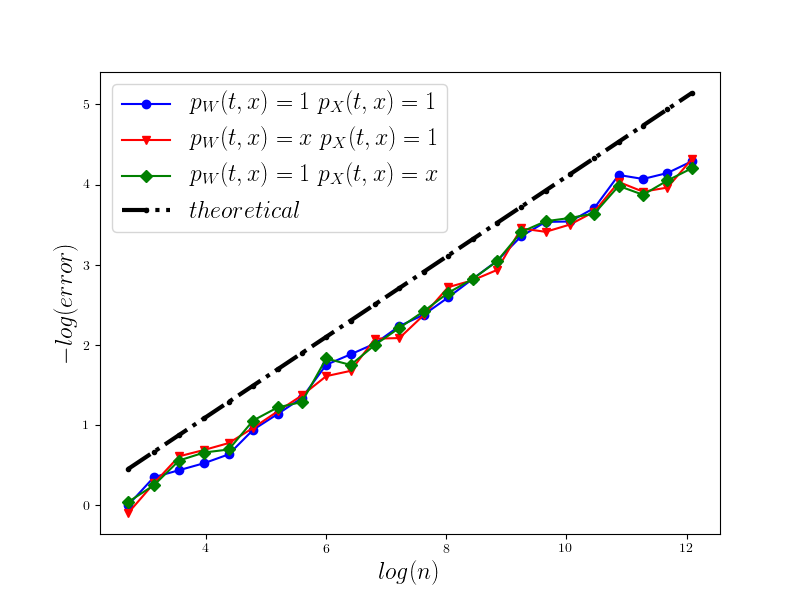}}

\caption{Error behavior for Riemann-Maruyama quadrature under exact/inexact information for problem \eqref{NR_PROB_1} - \eqref{NR_PROB_4} where $\delta_1, \delta_2$ are on the level $n^{-1/2}$.
}
\label{fig7}
\end{figure}

\begin{figure}
\includegraphics[width=0.8\linewidth]{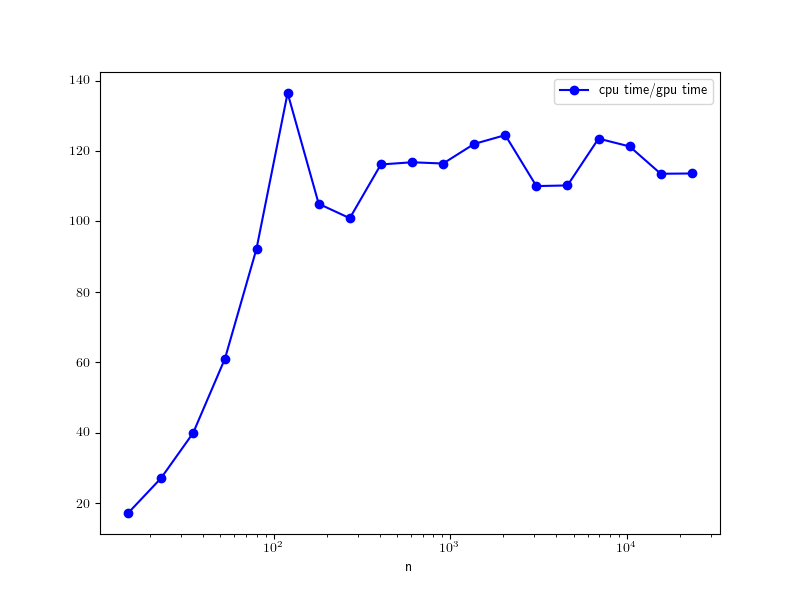} 
\caption{Performance speedup observed for GPU ({\em NVIDIA Tesla P100}) vs CPU ({\em Intel Broadwell})}
\label{fig5}
\end{figure}

\subsection{Problems}
For the test purposes we analyze following integration problem
\begin{equation}
	\mathcal{I}(X,W)=\int\limits_0^T X(t)\rd W(t)
\end{equation}
and we consider the following four examples
\begin{align}
	X_1(t) &= W(t),\qquad \mathcal{I}(X_1,W) = \frac{1}{2}W^2(T) - \frac{1}{2}T, \label{NR_PROB_1}\\
	X_2(t) & = W_2(t),\label{NR_PROB_2}\\
	X_3(t) &= f(S(t))S(t), \qquad f(x) = \max\{0, K-x\}, \qquad S(t)=S(0)e^{-\frac{1}{2}\sigma^2t+\sigma W(t)} \label{NR_PROB_3}\\
	\nonumber	& \mbox{where } K = 9, \sigma = 1  , S(0) = 1 ,\\
	X_4(t) & = N(t) e^{W(t)} \label{NR_PROB_4}, 
\end{align}
where $N=\{N(t)\}_{t\in [0,T]}$ is a Poisson process with insensitivity $\lambda= 5$ and $W_2=\{W_2(t)\}_{t\in [0,T]}$ is a standard one-dimensional Wiener process, both independent from $W$. We also apply $\mathcal{A}_n^{RM}$ to the weak approximation of the following scalar SDE
\begin{equation}
	\label{PROBLEM_SDE1}
		\left\{ \begin{array}{ll}
			dY(t)=\mu Y(t)\rd t+W_2(t)\rd W(t), &t\in [0,T], \\
			Y(0)=0, 
		\end{array}\right.
\end{equation}
where $\mu=3$.  The exact solution of (\ref{PROBLEM_SDE1}) leads to the quadrature problem, since
\begin{equation}
	Y(T)=\mathcal{I}(X,W)=\int\limits_0^T e^{\mu (T-t)}W_2(t)\rd W(t),
\end{equation}
where $\displaystyle{X(t)=e^{\mu (T-t)}W_2(t)}$, see \cite{esikr}, \cite{KP}. We use GPU implementation of the Riemann-Maruyama quadrature in order to compute an approximation of the following expectation
\begin{equation}
\label{weak_app_Y}
\mathbb{E}(f(Y(T)))=\mathbb{E}\Biggl(f\Bigl(\int\limits_0^T e^{\mu (T-t)}W_2(t)\rd W(t)\Bigr)\Biggr)
\end{equation}
 for $f=f(x)$ given as in (\ref{NR_PROB_3}) with $K=2$. Computation of (\ref{weak_app_Y}) corresponds to derivative pricing, where the price of the underlying risky asset is described by (\ref{PROBLEM_SDE1}).
 
 The approximation to $\displaystyle{\mathbb{E}(f(Y(T)))}$ is defined by
 \begin{equation}
	\frac{1}{M}\sum\limits_{j=1}^M f(\mathcal{A}_{n,j}^{RM}(\tilde X,\tilde W)),
 \end{equation}
 where $M$ is a number of independent copies of  $\mathcal{A}_{n}^{RM}(\tilde X,\tilde W)$. Due to the strong law of large numbers we get for all $n\in\mathbb{N}$
\begin{equation}
	\frac{1}{M}\sum\limits_{j=1}^M f(\mathcal{A}_{n,j}^{RM}(\tilde X,\tilde W))\to \mathbb{E}(\mathcal{A}_{n}^{RM}(\tilde X,\tilde W)), \quad \hbox{a.s.,}
\end{equation}
as $M\to +\infty$. Moreover, since $f:\mathbb{R}\to\mathbb{R}$ is a Lipschitz function and $X\in F_{L}^{1/2,2,q}$ with $q>2$, the standard arguments and Theorem \ref{main_thm} (i) yield the following estimate for averaged weak error, where  $(\tilde X, \tilde W)\in V_X(\delta_1)\times\mathcal{W}_s(\delta_2)$ and $\delta_1,\delta_2\in [0,1]$,
 \begin{align}
 \label{weak_err_rm}
	\Bigl\|\mathbb{E}(f(\mathcal{I}(X, W)))-\frac{1}{M}\sum\limits_{j=1}^M f(\mathcal{A}_{n,j}^{RM}(\tilde X,\tilde W))\Bigl\|_2&\leq C_1\cdot\|\mathcal{I}(X, W)-\mathcal{A}_n^{RM}(\tilde X,\tilde W)\|_2\notag\\
	&+2M^{-1/2}\|f(\mathcal{A}_n^{RM}(\tilde X,\tilde W))\|_2\notag\\
	&\leq C_2(n^{-1/2}+\delta_1+\delta_2+M^{-1/2}).
 \end{align}

\subsection{Noise}
For the purpose of testing we analyze following disturbing functions
\begin{align*}
	p_1(t,x) &= 1, \\
	p_2(t,x) &= x, \\
	p_3(t,x)  & = xt^2.
\end{align*}

It is worth to mention that the noise function $p_1$ corresponds to the standard absolute deterministic noise and $p_2$, $p_3$ are related to the standard relative error. The latter can be connected with the computation precision. There is a trend now observed in computations, e.g. for deep learning, where the computations are conducted in lower precision in order to gain huge computation speedup. The novel GPU architectures (e.g. NVIDIA Volta) are designed with some dedicated accelerators for single or half precision operations. 

For each test the information about the analyzed precision level $\delta_1$ and $\delta_2$ will be given. All the tests were conducted with $r=2$. 

\subsection{Error criterion}
For the problem \eqref{NR_PROB_1} we know the exact value of the solution, therefore we can have the following error estimate 

$$ \| \mathcal{I}(X_1,W)  - \mathcal{A}_n^{RM}(\tilde X_1, \tilde W) \|_2 \approx \left( \frac{1}{M} \sum_{j=1}^M \left( \mathcal{I}_j(X_1,W) - \mathcal{A}_{n,j}^{RM}(\tilde X_1, \tilde W)  \right)^2\right)^{1/2},
$$
where $M=2048$ corresponds to the number of computed independent realizations under given precision levels. 
In case of the problems (\ref{NR_PROB_2})-(\ref{PROBLEM_SDE1}), the exact solution is not known, hence in order to analyze the algorithm error, we need to compare the obtained result with the result obtained on the same trajectories for denser mesh. In our tests, as the expected convergence ratio is of no less than $0.5$, it is reasonable to have thousand times more points. That leads to following error estimation formula, used for (\ref{NR_PROB_2})-(\ref{NR_PROB_4})

\begin{equation}\label{NUM_ERR_INACC} 
	\| \mathcal{I}(X_m,W)  - \mathcal{A}_n^{RM}(\tilde X_m, \tilde W) \|_2 \approx \left( \frac{1}{M} \sum_{j=1}^M \left( \mathcal{A}_{L\cdot n,j}^{RM}(X_m, W) -\mathcal{A}_{n,j}^{RM}(\tilde X_m, \tilde W) \right)^2\right)^{1/2},
\end{equation}

where $m\in\{2,3,4\}$ and $L = 1000$.  
For (\ref{PROBLEM_SDE1}) we use the following quantity
 \begin{equation}
	\Bigl|\frac{1}{M}\sum\limits_{j=1}^M f(\mathcal{A}_{n,j}^{RM}(\tilde X,\tilde W))-\frac{1}{M}\sum\limits_{j=1}^M f(\mathcal{A}_{L\cdot n,j}^{RM}(X,W))\Bigl|
 \end{equation}

  as the approximation of the weak error
  \begin{equation}
	\Bigl|\mathbb{E}(f(\mathcal{A}_n^{RM}(\tilde X, \tilde W)))-\mathbb{E}(f(Y(T)))\Bigl|.
  \end{equation}
  
\subsection{Results}
In  Figure \ref{fig1} we present the behavior of the error for the Riemann-Maruyama quadrature $\mathcal{A}^{RM}_n$  for problem \eqref{NR_PROB_1}. The numerical results are compared with the theoretical rate of convergence obtained for the algorithm, i.e. we present the effect of changing the precision levels $\delta_1$ and $\delta_2$. In Figures \ref{fig2}-\ref{fig4} we present behavior of the error for the Riemann-Maruyama quadrature $\mathcal{A}^{RM}_n$ for problems \eqref{NR_PROB_2}-\eqref{NR_PROB_4}. (The  errors are measured accordingly to \eqref{NUM_ERR_INACC}.) From the Figure \ref{fig7} we see that if $\delta_1, \delta_2$ are on the level $n^{-1/2}$ then the  Riemann-Maruyama quadrature preserves the error $\mathcal{O}(n^{-1/2})$, known from the case when the information is exact. The results confirm the necessity of tending with the precision parameters to zero in order to maintain the convergence rate for the Riemann-Maruyama quadrature.

Results for the weak approximation are given at Figure \ref{fig6}.

\subsection{GPU implementation}
Below we present pseudo-code for the GPU implementation of the algorithm $\mathcal{A}_n^{RM}$. This algorithm is designed for the case where we wish to compute multiple realizations of $\mathcal{A}_n^{RM}$, returning the array of results. That algorithm, because of straightforward usage of parallel programming, enabled significant computational improvement of using graphics processing units. Moreover, additional speedup can be observed for using GPUs also for generating normally distributed numbers. Hence, it is suitable for e.g. derivative pricing. 


\begin{algorithm}[h]
   \caption{Riemann--Maruyama Quadrature}
    \begin{algorithmic}[1]
      \Function{$\mathcal{A}_n^{RM}$}{$n$, $T$, $M$, $\tilde X$, $\tilde W$} \Comment {$\tilde Z=(\tilde Z_1, \tilde Z_2, \ldots, \tilde Z_M)$, $\tilde Z\in\{\tilde X, \tilde W\}$}
        \State w = w\_{prev}= res = 0 \Comment{all variables are arrays of length M}
        \For{j = 1 to M} \Comment{this loop can be processed in parallel in independent threads}
        \State x[j] = $\tilde X_j(0)$
        \For{i = 0 to {n}}
        \State t = T*i/n		
            	\State w[j] = $\tilde W_j(t)$ 		
		\State res[j] $=$ res[j] + x[j]*(w[j]-w\_prev[j]) 
		\State x[j] = $\tilde X_j(t)$ 
		\State w\_{prev}[j] = w[j]   
        \EndFor 
        \EndFor\\
        \Return res
\EndFunction
\end{algorithmic}
\end{algorithm}


%

In our experiments we compared the performance of the algorithm $\mathcal{A}_n^{RM}$ for both GPU and CPU implementations. For GPU implementation we used 32 blocks and 64 threads for problems \eqref{NR_PROB_1}, \eqref{NR_PROB_2}, \eqref{NR_PROB_4} and 512 threads for problem \eqref{NR_PROB_3}. The CPU performance was tested using 8 threads. The used hardware was GPU -- {\em NVIDIA TESLA P100 (Maxwell)}, CPU -- {\em Intel Xeon E5-2680v4 (Broadwell)}. The speedup comparison for problem (\ref{NR_PROB_3}) is presented in Figure \ref{fig5}. As we can see it is possible to have speedup of level 100x.


\section{Conclusions} \label{sec:con}
We investigated the problem of efficient approximation of It\^o integrals under inexact information about the Wiener process and an integrand. We showed that for certain precisions ($\delta_1=\delta_2\geq 0$) the Riemann-Maryama quadrature rule is optimal. We also proposed  GPU implementation of the algorithm that is suitable for practical purposes.
\newline\newline
{\bf Acknowledgments.}

The authors were partially supported by the Faculty of Applied Mathematics AGH UST dean grant for PhD students and young researchers within subsidy of Ministry of Science and Higher Education, with grant numbers as follows: 15.11.420.038/18 (Andrzej Ka\l u\.za), 15.11.420.038/2 (Pawe\l \ M. Morkisz), and 15.11.420.038/1 (Pawe\l \ Przyby\l owicz).


\begin{thebibliography}{22}
\bibitem{delmey}
C. Dellacherie, P-A. Meyer, Probabilities and Potential, Hermann, 1978.

\bibitem{esikr}
M. Eisenmann, R. Kruse, Two quadrature rules for stochastic It\^o integrals with fractional Sobolev regularity, https://arxiv.org/abs/1712.08152.

\bibitem{Fatahalian}
K. Fatahalian, J. Sugerman,  P. Hanrahan, Understanding the efficiency of GPU algorithms for matrix-matrix multiplication, Graphics Hardware (2005), 133--137.

\bibitem{alex}
N. Whitehead, A. Fit-Florea, Precision \& Performance: Floating Point and IEEE 754 Compliance for NVIDIA GPUs, NVIDIA, 2011.

\bibitem{legall}
J-F. Le Gall, {\em Brownian Motion, Martingales, and Stochastic Calculus}, Springer, 2016.
 
\bibitem{hein1}
S. Heinrich, Lower complexity bounds for parametric stochastic It\^o integration, Preprint, https://www.uni-kl.de/AG-Heinrich/papers/lowpsint17.pdf,  2017. 

\bibitem{heindaun}
S. Heinrich, T. Daun, Complexity of Banach space valued and parametric stochastic It\^o integration, {\em J. Complex.} {\bf 40}  (2017), 100--122.
 
\bibitem{hert}
P. Hertling,  Nonlinear Lebesgue and It\^o integration problems of high complexity, {\em J. Complex.} {\bf 17}  (2001), 366--387.

 
\bibitem{JoC4}
B. Kacewicz, M. Milanese, A. Vicino, Conditionally optimal algorithms and estimation of reduced order models,
{\em J. Complexity} {\bf 4} (1988), 73--85.


\bibitem{KaPl90}
B. Kacewicz, L. Plaskota, On the minimal cost of approximating linear problems based on information with deterministic noise, {\em Numer. Funct. Anal. and Optimiz.} {\bf 11} (1990), 511-528.

\bibitem{KaPr16}
B. Kacewicz, P. Przyby{\l}owicz, On the optimal robust solution of IVPs with noisy information, {\em Numer. Algor.} {\bf 71} (2016),  505--518.

\bibitem{KARSHR}I. Karatzas, S. E. Shreve, \textit{Brownian Motion and Stochastic Calculus}, 2n edition, Springer - Verlag, New York, 1991.

\bibitem{KP} P. E. Kloeden, E. Platen, \textit{Numerical Solution of Stochastic Differential Equations}, Springer-Verlag Berlin, Heidelber, third edition, 1999.

\bibitem{Kruger}
J. Kr\"{u}ger, R. Westermann, R\"{u}diger, Linear Algebra Operators for GPU Implementation of Numerical Algorithms, {\em ACM Trans. Graph.} {\bf 22} (2003), 908--916.

\bibitem{Langdon} W. B.  Langdon,  Brian Y. H. Lam, J. Petke, M.  Harman. Improving CUDA DNA Analysis Software with Genetic Programming. Proceedings of the 2015 on Genetic and Evolutionary Computation Conference - GECCO '15. (2015), 1063–1070.

\bibitem{NV_Volta} NVIDIA, \textit{NVIDIA Tesla V100 GPU Architecture}, available online http://images.nvidia.com/content/volta-architecture/pdf/volta-architecture-whitepaper.pdf.




\bibitem{MoPl16}
 P. M. Morkisz, L. Plaskota, Approximation of piecewise H\"older functions from inexact information, {\em J. Complex.} {\bf 32} (2016), 122--136.
 
\bibitem{MoPr17} 
P. M. Morkisz, P. Przyby\l owicz, Optimal pointwise approximation of SDE's from inexact information, {\em J. Comp. and Appl. Math.} {\bf 324} (2017), 85--100.


\bibitem{Pla96}
  L. Plaskota, {\em Noisy Information and Computational Complexity}, 
   Cambridge Univ. Press, Cambridge, 1996.

\bibitem{Pla14}
  L. Plaskota, Noisy information: optimality, complexity, tractability, 
  in {\em Monte Carlo and quasi-Monte Carlo Methods 2012}, 
  J. Dick, F.Y. Kuo, G.W. Peters, I.H. Sloan (Eds.), Springer 2013,  173--209. 

\bibitem{PP1}
P. Przyby\l owicz, Linear information for approximation of the It\^o integrals. 
{\it Numer. Algorithms} {\bf 52} (2009), 677--699. 


\bibitem{PP2} 
P. Przyby\l owicz, Adaptive It\^o--Taylor algorithm can optimally approximate the It\^o integrals 
of singular functions. {\it J. Comp. Appl. Math.} {\bf 235} (2010), 203--217.

\bibitem{PP3} 
P. Przyby\l owicz, Optimal sampling design for approximation of stochastic It\^o integrals with 
application to the nonlinear Lebesgue integration. {\it J. Comp. Appl. Math.} {\bf 245} (2013), 10--29.

\bibitem{Riesinger}
C. Riesinger, T. Neckel, F. Rupp, Solving random ordinary differential equations on GPU clusters using multiple levels of parallelism, {\it SIAM J. Sci. Comput.} {\bf 38(4)} (2016), pp. C372--C402.

\bibitem{Riesinger2}
C. Riesinger, T. Neckel, F. Rupp, A. P. Hinojosa, H.-J. Bungartz, GPU optimization of pseudo random number generators for random ordinary differential equations, {\it Procedia Computer Sci} {\bf 29} (2014), pp. 172-183.




\bibitem{Ryoo}
S. Ryoo, C. I. Rodrigues, S. S. Baghsorkhi, S. S. Stone, D. B. Kirk, W. W. Hwu, Optimization principles and application performance evaluation of a multithreaded GPU using CUDA, Proceedings of the 13th ACM SIGPLAN PPoPP, (2008), 73-82.






\bibitem{TWW88}
  J.F. Traub, G.W. Wasilkowski,  H. Wo\'zniakowski,
  {\em Information-Based Complexity}, Academic Press, New York, 1988.

\bibitem{waswo}
G. W. Wasilkowski, H.  Wo\'zniakowski, On the complexity of stochastic integration, {\em Math. Comp.} {\bf 70} (2001), 685-698.


\bibitem{Wer96}
  A.G. Werschulz, The complexity of definite elliptic problems with noisy data. {\em J. Complex.} {\bf 12} (1996), 440-473.
  
\bibitem{Wer97}
  A.G. Werschulz, The complexity of indefinite elliptic problems with noisy data. {\em J. Complex.} {\bf 13} (1997), 457-479.
\end{thebibliography}
\end{document}